\documentclass[10pt,a4paper]{amsart}
\usepackage{amsmath, amssymb, amsthm}
\usepackage{mathtools}
\usepackage{amsaddr}
\usepackage{geometry}
\usepackage[utf8]{inputenc}
\usepackage{lmodern}
\usepackage[T1]{fontenc}
\usepackage{enumitem}
\usepackage{graphicx}
\graphicspath{{./figs/}}
\usepackage{caption,subcaption}
\usepackage{colonequals}
\usepackage{siunitx}
\usepackage{tikz}
\usetikzlibrary{arrows.meta}

\usepackage{pgfplots}
\pgfplotsset{compat=newest}

\newcommand{\ordernum}[1]{\num[round-mode=places,
                               round-precision=3,
                               round-pad=false,
                               exponent-mode=fixed,
                               fixed-exponent=0]{#1}}

\newcommand{\timenum}[1]{\num[round-mode=places,
                              round-precision=3,
                              round-pad=false,
                              exponent-mode=fixed,
                              fixed-exponent=0]{#1}}

\usepackage[linesnumbered,ruled,vlined]{algorithm2e}

\usepackage{longtable} 
\usepackage{ltablex}
\usepackage{xltabular}
\usepackage{array, ltablex, caption, multirow, makecell}

\def\R{\mathbb{R}}
\def\cA{\mathcal{A}}
\def\cI{\mathcal{I}}
\def\cL{\mathcal{L}}
\def\cO{\mathcal{O}}
\def\cP{\mathcal{P}}
\def\cS{\mathcal{S}}
\def\bn{\mathbf{n}}
\def\bu{\mathbf{u}}
\def\bv{\mathbf{v}}
\def\bphi{\boldsymbol{\phi}}
\def\bvarphi{\boldsymbol{\varphi}}

\def\Th{\mathcal{T}_h}
\def\Ih{\mathcal{I}_h}
\def\DD{{\rm D}}
\def\uD{\bu_\DD}
\newcommand{\dv}[1]{\,{\mathrm d}#1}

\theoremstyle{definition}
\newtheorem{definition}{Definition}
\newtheorem{problem}{Problem}
\newtheorem{remark}[definition]{Remark}
\theoremstyle{plain}
\newtheorem{lemma}[definition]{Lemma}
\newtheorem{theorem}[definition]{Theorem}

\makeatletter
\@namedef{subjclassname@2020}{\textup{2020} Mathematics Subject Classification}
\makeatother

\usepackage{todonotes}

\usepackage{hyperref}

\usepackage[sort,nocompress]{cite}

\DeclarePairedDelimiter{\abs}{\lvert}{\rvert}
\DeclarePairedDelimiter{\norm}{\lVert}{\rVert}

\def\Exp{\operatorname{Exp}}

\def\CC{{C\nolinebreak[4]\hspace{-.05em}\raisebox{.4ex}{\tiny\bf ++}}}

\usepackage[most]{tcolorbox}
\definecolor{lightgrey}{rgb}{0.9, 0.9, 0.9}
\definecolor{lightblue}{rgb}{20,140,190}
\definecolor{cornflowerblue}{rgb}{0.39, 0.58, 0.93}
\newtcolorbox{problembox}{colback=lightgrey, colframe=black}

\begin{document}
\title{Benchmarking Numerical Algorithms for Harmonic Maps into the Sphere}
\author{Sören Bartels$^1$, Klaus Böhnlein$^2$, Christian Palus$^1$, Oliver Sander$^2$}
\address{$^1$ Department of Applied Mathematics, University of Freiburg \\
  $^2$ Institute of Numerical Mathematics, Technische Universität Dresden}
\date{\today}
\begin{abstract}
We numerically benchmark methods for computing harmonic maps into the unit sphere,
with particular focus on harmonic maps with singularities.
For the discretization we compare two different approaches,
both based on Lagrange finite elements. While the first method enforces the
unit-length constraint only at the Lagrange nodes, the other one adds a pointwise projection
to fulfill the constraint everywhere.
For the solution of the resulting algebraic problems we compare a nonconforming
gradient flow with a Riemannian trust-region method. Both are energy-decreasing
and can be shown to converge globally to stationary points of the discretized Dirichlet energy.
We observe that while the
nonconforming and the conforming discretizations both show similar behavior
for smooth problems, the nonconforming discretization handles singularities better.
On the solver side,
the second-order trust-region method converges after few steps, whereas the number
of gradient-flow steps increases proportionally to the inverse grid element diameter.
\end{abstract}

\maketitle

\section{Introduction}

Harmonic maps are stationary configurations of the Dirichlet energy~\cite{helein2008harmonic,eells_lemaire:1995}.
In this work we focus on maps into the unit sphere
$S^{m-1} \colonequals \{ x \in \R^m \; : \; \abs{x} = 1\}$.  More formally,
given a Lipschitz domain $\Omega \subset \R^n$ with $n \in \{2,3\}$,
we seek stationary configurations $\bu \colon \Omega \to S^{m-1}$ with $m\in\{2,3\}$, $m\geq n$, of the Dirichlet energy
\begin{equation}\label{eq:dirichlet-energy}
  E[\bu] \colonequals \frac{1}{2} \int_\Omega \abs{\nabla \bu}^2 \dv{x}.
\end{equation}
Equivalently, we will sometimes regard the problem as looking for vector fields
$\bu : \Omega \to \R^m$ that make $E[\cdot]$ stationary while fulfilling the constraint
\begin{equation}\label{eq:unit-constr}
  |\bu|^2 = 1 \text{ almost everywhere in } \Omega.
\end{equation}
For well-posedness we require Dirichlet boundary conditions
\begin{equation}
  \bu = \uD \qquad \text{on $\partial \Omega$}, \label{eq:boundary-cond}
\end{equation}
for a function $\uD : \Omega \to S^{m-1}$ of suitable smoothness~\cite[Chap.\,3.4]{helein2008harmonic}.

Harmonic maps into spheres have numerous practical applications.
They appear in models of liquid crystal materials, which, on a
microscopic level, consist of rod-shaped molecules that exhibit a natural desire for mutual alignment.
Popular macroscopic liquid crystal models are based on the Oseen--Frank energy,
which is the Dirichlet energy of the field of molecule orientations.
The emergence of new technologies in the manufacturing of liquid crystal compound materials
as well as new ideas for technical applications~\cite{KohChe13,deHaan_et_al14,WhiBro15,Ware_et_al16}
have led to an ongoing interest also in corresponding simulation schemes~\cite{NocWalZha17, Walker20, BarGriNeuParPal23}.

Harmonic maps into spheres also play a role in micromagnetics, which models
orientation fields of magnetic dipoles. The Dirichlet energy is the simple-most representative
of a family of different energy-based models~\cite{desimone2004recent}.
It also serves as a prototype energy that already captures a range of interesting effects.
Furthermore, harmonic maps serve as the basis
for more complicated energies such as the Ginzburg--Landau model or the chiral skyrmions discussed in~\cite{melcher:2014}.
In this context singular harmonic maps are particularly interesting, since they arise as limits of minimizers of the Ginzburg--Landau energy~\cite{monteil2021ginzburg}.

Finally, from a mathematical point of view harmonic maps are interesting in their own right,
and a considerable body of literature exists with investigations of mathematical properties
of harmonic maps. Overviews and further literature can be found in~\cite{helein2008harmonic,struwe2000variational,BreCorLie86}.

\medskip

The construction of finite element (FE) methods for the approximation of harmonic maps
requires special care. The central challenge is the handling of the non-Euclidean image space $S^{m-1}$,
or, equivalently, of the nonlinear, nonconvex constraint $\abs{\bu(x)}^2 = 1$ for
(almost) all $x$ in the domain $\Omega$.  This affects
both the discretization of the problem as well as solution strategies for the resulting 
algebraic systems.
In the last decades various approaches have been proposed that can be employed in the numerical approximation of harmonic maps.
For example, several authors have proposed finite-difference approximations of the Dirichlet energy
\cite{cohen1987minimum,weinmann2014total}, as well as point-relaxation methods~\cite{lin1987numerical}
and gradient-type methods~\cite{Alouges97,bartels2005robust}.
For the constraint, approaches like parametrizations~\cite{vese2002numerical}, Lagrange multipliers~\cite{cohen1987minimum,huWinther2009saddle} and penalization~\cite{liu2000approximation,du1992analysis} have been proposed.
Several works also treat numerical methods for more general problems such as
$p$-harmonic maps \cite{vese2002numerical}
or fractional harmonic maps~\cite{antil2023approximation}.

While all these methods warrant thorough benchmarking, for reasons of space
we limit this paper to two particular discretizations and two particular
solver algorithms.
In~\cite{Bartels05,Bartels16}, Bartels and coworkers proposed a nonconforming discretization method that consists
of $\R^m$-valued Lagrange finite elements that are constrained to fulfill the
unit-length constraint only at the Lagrange points.
Stationary points of the discretized energy~\eqref{eq:dirichlet-energy} are approximated
via an implicit discrete gradient flow employing a linearization of the unit-length constraint~\eqref{eq:unit-constr}.
This solver approach is nonconforming in the sense that the iterates slowly accumulate a violation
of the constraints. Bartels et al.\ showed, however, that this constraint violation remains bounded in terms of the time step size $\tau$,
and that the discrete solutions weakly converge to stationary points of the Dirichlet energy
as the grid element size $h$ and the time step size~$\tau$ go to zero.
The energy monotonicity property of the discrete gradient flow
implies a convergence rate of $\mathcal{O}\big((\tau k)^{-\frac{1}{2}}\big)$ for the norms of the corrections, where $k$ denotes the number of solver steps.

The second approach has been proposed by Sander et al.~\cite{Sander12,Sander16},
who aimed at a completely
conforming method. The authors construct so-called geometric finite elements,
which are generalizations of piecewise polynomial functions (of arbitrary order) that map
into $S^{m-1}$ at any point in the domain.
For problems with sufficiently smooth solutions,
Hardering et al.\ showed optimal convergence rates for the $H^1$- and $L^2$-discretization errors
for these elements, for any polynomial order~\cite{GFE, GFEprojection}.
Remarkably, such finite elements can also represent singular maps to a certain extent,
which suggests to employ them for simulating harmonic maps with singularities.
To solve the discrete problems Sander et al.\ interpreted them
as optimization problems on a product manifold, and solved them using a Riemannian
trust-region method~\cite{absil2009optimization,Sander12}. The convergence of this algebraic solver follows from general results for
optimization methods on manifolds~\cite{absil2009optimization}.

The goal of this paper is to compare the practical properties of the two different 
discretizations and solver algorithms for the numerical approximation of harmonic maps into the unit sphere.
The overall outcome is not a priori clear: While the conforming methods preserve more of the
mathematical structure of the problem, the nonconforming methods are simpler to implement.
Also, the behavior in the presence of singularities is hard to predict; indeed,
previous works already observe that computing the exact placement of singularities
is difficult~\cite{huWinther2009saddle}.
As both solver methods are based on energy minimization they are unlikely to find anything but locally minimizing harmonic maps.

To test the relative merits of the methods we define a set of benchmark problems.
This set includes smooth harmonic maps, but also maps with different types of singularities.
In our measurements we focus on the discretization error, the constraint violation of the
nonconforming methods, and the solver speed.
We find that the approximation power of the two discretizations is the same
for smooth problems, but the conforming discretization shows slightly better
convergence orders than the nonconforming one in the presence of singularities.
For problems with smooth solutions, discretization error convergence orders
are clearly observed, and they match the theoretical predictions.
For problems with singular solutions, it is much more difficult
to determine an order, and the precise behavior depends on the exact position
of the singularity. Both solvers converge, but for three-dimensional problems with a singularity,
we observe that the placement of the final singularity can depend on the initial iterate.
This effect, which was already noted in~\cite{huWinther2009saddle}
for two-dimensional problems, is more pronounced
for the conforming discretization than for the nonconforming one,
and remains to be addressed in future work.

As predicted by theory, the constraint violation produced by the gradient-flow solver
remains bounded as long as the step size is chosen to be proportional
to the element diameter. The values we measured are small enough to be
unproblematic for most practical applications. By construction, the trust-region solver
does not introduce any constraint violation at all.

Considerable differences show when comparing the solver speeds. The iteration numbers
of the trust-region method are low, and in most situations they seem to remain bounded
independently of the grid resolution.  The gradient-flow solver, on the other hand,
has to tie its step size to the grid element diameter. That way, the total number
of iterations increases as the grid is refined. As the iterations of both methods
have comparable cost, this leads to a large speed advantage for the trust-region
solver.
Experiments that combine the gradient-flow method with a Newton solver
can be found in~\cite{bartels:2009}.

\medskip

The paper is structured as follows: Chapter~\ref{sec:characterization} briefly
reviews the different notions of harmonic maps used in this text.
Chapter~\ref{sec:discretizations} then introduces the two discretization
methods, and Chapter~\ref{sec:solvers} does the same for the algebraic solvers.
The different benchmark problems are presented in Chapter~\ref{sec:benchmarks-problems}.
Finally, the last two chapters contain the actual numerical results,
with Chapter~\ref{sec:discretization_benchmarks} investigating the discretizations, and Chapter~\ref{sec:solver-benchmark}
the solvers.

\section{Harmonic maps into the sphere}\label{sec:characterization}

Let $\Omega$ be an open, bounded domain in $\R^n$.
We use the standard notation
$H^{k}(\Omega; \R^m)$ for vector-valued Sobolev spaces defined on $\Omega$ with  $k\geq 1$,
and we denote the $k$-th order Sobolev vector fields with vanishing boundary trace by $H^{k}_0(\Omega; \R^m)$.
For sphere-valued problems, we introduce the subspace
\begin{align*}
H^k(\Omega;S^{m-1}) &\colonequals \Big\{ \bv\in H^k(\Omega;\R^m) : \bv(x) \in S^{m-1} \text{ a.e.} \Big\},
\end{align*}
and its subspace of functions that satisfy the boundary conditions~\eqref{eq:boundary-cond} in a trace sense
\begin{align*}
H^k_D(\Omega;S^{m-1} )
&\colonequals
\Big\{ \bu \in H^k(\Omega;S^{m-1} ) \; : \; \bu|_{\partial \Omega} = \uD|_{\partial \Omega} \Big\}.
\end{align*}

There are various related definitions of harmonic maps.
This text uses three of them, which we review here briefly.  More details
can be found in~\cite{helein2008harmonic,BreCorLie86,eells_lemaire:1995}.
In the following, $(\cdot,\cdot)$ denotes the standard $L^2$ scalar product.

\begin{definition}[Stationary harmonic map]
\label{def:stationary_harmonic_map}
A function $\bu \in H^1(\Omega; S^{m-1})$
is called \emph{stationary harmonic map}
if it satisfies the weak Euler--Lagrange equation
\begin{equation}\label{eq:geometric-euler-lagrange}
(\nabla \bu, \nabla \bvarphi) = 0
\end{equation}
for all $\bvarphi \in H^1_0(\Omega; \R^m)$ with $\bu \cdot \bvarphi = 0$ a.e. in $\Omega$.
\end{definition}

The tangentiality condition $\bu \cdot \bvarphi = 0$ can be avoided by the equivalent formulation
\begin{equation*}
(\nabla \bu, \nabla \bvarphi) = (|\nabla \bu|^2 \bu, \bvarphi)
\end{equation*}
where we now test with all $\bvarphi \in H^1_0(\Omega; \R^m) \cap L^\infty(\Omega; \R^m)$, see~\cite{Bartels15_book}.

\begin{definition}[(Locally) minimizing harmonic map]
A map $\bu \in H^1(\Omega;S^{m-1})$
 is called \emph{(locally) minimizing harmonic map} if it is a (local) minimizer of
the Dirichlet energy~\eqref{eq:dirichlet-energy}.
\end{definition}

In general, a minimizing harmonic map of sufficient smoothness is also
a stationary harmonic map, but the converse is not always the case.

The characterization~\eqref{eq:geometric-euler-lagrange} can be generalized via partial integration to cover cases
of lesser regularity.

\begin{definition}[Distributional harmonic map]
A map $\bu \in L^\infty(\Omega; \R^m)$ with $\abs{\bu}=1$ a.e.\ in~$\Omega$ is called
\emph{distributional harmonic map} if
\begin{equation*}
  (\bu, \Delta \bvarphi) = 0
\end{equation*}
for all $\bvarphi \in C_{0}^\infty(\Omega; \R^m)$ with $\bu \cdot \bvarphi = 0$ almost everywhere in $\Omega$.
\end{definition}

It is a particularity of sphere-valued problems
that the set of continuous maps $\Omega \to S^{m-1}$ is not connected.
In other words, given two continuous functions $\bv_1, \bv_2 : \Omega \to S^{m-1}$
with identical Dirichlet boundary values
it may not be possible to continuously deform one into the other.
The underlying reason is the fact that the homotopy groups of the sphere are not all trivial.
Indeed, maps from a simply connected $n$-dimensional domain $\Omega \subset \R^n$ into $S^{m-1}$ are
closely related to the $n$-th homotopy group of $S^{m-1}$. This group is trivial
in the case $n = m = 2$, but isomorphic
to~$\mathbb{Z}$ in the cases $n = m = 3$ and $n = 2, m = 3$ considered here.
By this isomorphism the connected components of the function spaces can be labeled by an integer
(which is is sometimes called topological degree or quantum number~\cite{thouless:1998}).
The disconnected nature of the function space is a
structure that may or may not be preserved by discretizations.

The issue of connectedness is more subtle for the Sobolev spaces $H^1(\Omega;S^{m-1})$
that underly finite element theory.  In particular, the homotopy classes
of $W^{1,p}$ may not all be weakly closed, and therefore some of them may not contain
minimizing harmonic maps~\cite[Chap.\,5.3]{helein2008harmonic}.
A discussion about the connected components
of $W^{1,p}(M,N)$ for general manifolds $M$, $N$ can be found in~\cite[Chap.\,3.3]{helein2008harmonic}.

\section{Discretizations} \label{sec:discretizations}

We now present
the two discretization approaches of~\cite{Bartels05,Bartels16} and~\cite{Sander12,Sander16} for
maps from the flat
domain $\Omega$ into $S^{m-1} \subset \R^m$.
The central challenge is the fact that it is impossible to satisfy the unit-length constraint~\eqref{eq:unit-constr}
everywhere in $\Omega$ if piecewise polynomials are used for the approximation.
In the following we assume $\Omega$ to be polyhedral, and we let $\Th$ be a shape-regular triangulation of $\Omega$ with triangles or tetrahedra
of diameter no larger than $h > 0$.
The set of polynomials with degree at most $p \ge 0$ on a simplex $T \in \Th$ is denoted with $P_p(T)$, and the space of continuous, piecewise polynomial $m$-vector fields is
\begin{equation*}
\cS^p(\Th; \R^m)
\colonequals
\Big\{ \bv_h \in C(\Omega; \R^m) \; : \; \bv_h|_T \in P_p(T)^m \text{ for all $T \in \Th$} \Big\}.
\end{equation*}
We use the notation $\Ih^p : H^2(\Omega;\R^m)\to\cS^p(\Th;\R^m)$ for the $p$-th order Lagrange
interpolation operator associated to a set of Lagrange points $\mathcal{L}_h$, of size $N \colonequals \abs{\cL_h}$.

\subsection{Geometrically nonconforming Lagrange finite elements} \label{sec:discret-nonconf}
A natural approach
to cope with the nonlinear nature of the image space $S^{m-1}$
is to approximate sphere-valued maps with functions from $\cS^p(\Th; \R^m)$
and require the unit-length constraint~\eqref{eq:unit-constr} only at the Lagrange points.
For first-order finite elements this has been proposed and analyzed in~\cite{Bartels05}.
We obtain the following admissible set, which also  enforces
the boundary conditions:
\begin{align*}
  \cA_{h}^\text{nc} \colonequals \Big\{
    \bv_h \in \cS^p(\Th; \R^m) : |\bv_h(z)| = 1 \text{ for all $z \in \mathcal{L}_h$ and $\bv_h(z) = \bu_\DD(z)$ for all $z \in \mathcal{L}_h \cap \partial \Omega$}
  \Big\}.
\end{align*}
This is a smooth nonlinear submanifold of the vector space $\cS^p(\Th; \R^m)$,
but observe that $\cA_h^\text{nc}$ is a connected set.
It therefore waives one of the central features of the actual solution space,
which consists of disconnected homotopy classes.

Discrete approximations of minimizing harmonic maps can now be defined as minimizers
of the Dirichlet energy~\eqref{eq:dirichlet-energy} in the discrete admissible set $\cA_h^\text{nc}$.
This works because while $\cA_h^\text{nc}$ is not a subset of the actual solution space $H^1(\Omega;S^{m-1})$,
the Dirichlet functional~\eqref{eq:dirichlet-energy} extends naturally to $H^1(\Omega;\R^m) \supset \cA_h^\text{nc}$.

The Definition~\ref{def:stationary_harmonic_map} of weak stationary harmonic maps
involves test functions.
As the space of admissible functions is nonlinear, the corresponding test function spaces differ
from point to point. Formally, the test functions for a map $\bv\in H^1(\Omega;S^{m-1})$
are elements of the tangent space $T_{\bv} H^1(\Omega;S^{m-1})$, and likewise for discrete maps $\bv_h \in \cA_h^\text{nc}$.
More practically, test functions are constructed as variations of admissible functions \cite{sander2016test}.
For a given $\bv_h \in \cA_h^\text{nc}$ let therefore $\gamma :[-\epsilon,\epsilon] \to \cA_h^\text{nc}$ be a
differentiable path with $\gamma(0) = \bv_h$. The set of test functions for $\cA_h^\text{nc}$
at~$\bv_h$ is the set of all functions $\bvarphi_h :\Omega \to \R^m$ that can be represented as
\begin{align*}
  \bvarphi_h = \frac{d\gamma(t)}{dt}\bigg|_{t=0}
\end{align*}
with such a path $\gamma$.
For the geometrically nonconforming discretization,
this construction results in a space of piecewise polynomial vector fields that are orthogonal
to $\bv_h$ at the Lagrange points, and zero at the boundary
\begin{align*}
  T_h^\textup{nc}(\bv_h)
  \colonequals \Big\{
    \bvarphi_h \in \cS^p(\Th; \R^m) : \bv_h(z)\cdot \bvarphi_h(z)  = 0 \text{ for all $z \in \mathcal{L}_h$},
    \; \bvarphi_h(z) = 0 \text{ for all $z \in \mathcal{L}_h \cap \partial \Omega$}
  \Big\}.
\end{align*}
Finding discrete stationary harmonic maps then means finding functions $\bu_h \in \cA_h^\text{nc}$
such that
\begin{equation*}
  (\nabla \bu_h, \nabla \bvarphi_h) = 0
\end{equation*}
for all $\bvarphi_h \in T_h^\textup{nc}(\bu_h)$. This corresponds to the continuous problem~\eqref{eq:geometric-euler-lagrange}.

\medskip

The nonconforming discretization of harmonic maps allows for different convergence theories.
All results in the literature concern the first-order case $p=1$ only, but the convergence theories
from~\cite{Bartels16} apply verbatim to higher-order discretizations.
Most generally, one can establish the $\Gamma$-convergence of the discretizations
of the Dirichlet energy~$E$ from~\eqref{eq:dirichlet-energy} defined by
\[
 E_h[\bu_h]
 \colonequals
 \begin{cases}
  E[\bu_h]  &\text{if $ \bu_h \in \cA_h^\text{nc}$},\\
  \infty    &\text{else},
\end{cases}
\]
to the continuous minimization problem.
A discrete compactness result then implies the accumulation of approximations at solutions of the continuous problem.
A proof for the case $m=3$ is given in \cite{Bartels15_book}.

\begin{theorem}[Discrete compactness,~{\cite[Theorem~7.6]{Bartels15_book}}]\label{thm:discrete-compactness}
  Let $p=1$ and let $(\bu_h)_{h>0} \subset H^1(\Omega; \R^3)$ be a bounded sequence such that each $\bu_h \in \mathcal{A}_h^\textup{nc}$ and
  \begin{equation*}
    (\nabla \bu_h, \nabla \bvarphi_h) = \mathcal{R}_h(\bvarphi_h)
  \end{equation*}
  for all $\bvarphi_h \in T_h^\textup{nc}(\bu_h)$ for functionals $\mathcal{R}_h \in H^1_0(\Omega; \R^3)'$ with $\mathcal{R}_h \to 0$ in $H^1_0(\Omega; \R^3)'$ as $h \to 0$.
  Then every weak accumulation point of $(\bu_h)_{h>0}$ is a stationary harmonic map.
\end{theorem}
\begin{proof}[Sketch of the proof]
  The satisfaction of the unit-length constraint by the accumulation points follows from standard interpolation and inverse estimates.
  The fact that weak limits are harmonic maps follows in the case $m=3$ from
  the equivalent characterization
  \begin{equation*}
    (\nabla \bu, \nabla[\bu \times \bphi]) = 0
  \end{equation*}
  for all $\bphi \in H^1_0(\Omega;\R^3) \cap L^\infty(\Omega;\R^3)$, by choosing
  the interpolant of $\bu_h \times \bphi$ for given $\bphi \in C^\infty(\overline{\Omega};\R^3)$ as a test function and then passing to the limit.
\end{proof}

This result can be extended to manifolds that are more general than the unit-sphere~\cite{Bartels10}.

\begin{remark}
In the proof of~Theorem~\ref{thm:discrete-compactness} one establishes the fact that for a sequence of discrete vector fields whose nodal values
belong to $S^{m-1}$, any weak accumulation point in $H^1(\Omega)$ has values in $S^{m-1}$ almost everywhere.
The assumption of the nodal values belonging to $S^{m-1}$ can be further weakened: Given a sequence of discrete vector fields whose nodal values
are not necessarily on $S^{m-1}$ but approach unit length as $h \to 0$,
it still follows that any accumulation point in $H^1(\Omega;\R^m)$ satisfies the constraint exactly; cf.~\cite{Bartels16}.
This allows to prove convergence of solutions obtained with
the nonconforming discrete gradient flow discussed in Section~\ref{sec:gradient-flow},
which does not preserve the constraint exactly.
\end{remark}

Recently, discretization error bounds have also been proved.
For piecewise linear Lagrange elements and $n=m=2$, optimal error estimates in the energy norm
have been derived for a corresponding saddle-point formulation in~\cite{huWinther2009saddle}.
This result is extended in~\cite{BaPaWa2023} to the cases $m=3$ and $n>2$.

\subsection{Geometrically conforming projection-based finite elements}
\label{sec:discret-conf}

The second discretization constructs
finite elements that map into the sphere $S^{m-1}$ everywhere, not just at the Lagrange points.
It does so by adapting the notion of a polynomial.
The construction has been proposed in~\cite{GFEprojection} and \cite{gawlikLeok2018embedding} under the name of
\emph{projection-based finite elements}, and is a member of
the larger family of \emph{geometric finite elements} \cite{GFE}.

In essence, projection-based finite elements are defined
by projecting Lagrange finite elements with nodal values in $S^{m-1}$ pointwise onto $S^{m-1}$. More formally, consider
the closest-point projection from $\R^m$ onto $S^{m-1}$ defined by
\begin{align*}
P:\R^m\setminus\{0\} \to S^{m-1},
\qquad
P(\xi)\colonequals \frac{\xi}{|\xi|}.
\end{align*}
This projection induces a superposition operator \cite{appell1990nonlinear}
\begin{align*}
  \cP : C(\Omega;\R^m) \to W^{k,q}(\Omega;S^{m-1})
  \qquad
  \cP \bv(x) \colonequals P(\bv(x)),
  \qquad
  \forall x \in \Omega,
\end{align*}
for suitable $k$ and $q$.
With its help we define the space of $p$-th order projection-based finite elements as
\begin{align*}
  \cS^{p,\textup{proj}}(\Th;S^{m-1})
  & \colonequals
  \Big\{ \cP\bv_h \; :\; \bv_h \in \cS^p(\mathcal T_h;\R^m)
    \text{ and $\bv_h(z) \in S^{m-1}$ for all $z \in \mathcal L_h$} \Big\}.
\end{align*}
A canonical interpolation operator into $\cS^{p,\textup{proj}}(\Omega;S^{m-1})$  is given by
\begin{align*}
  \Ih^{p,\textup{proj}} : C(\Omega;\R^m)&\to\cS^{p,\textup{proj}}(\Th;S^{m-1}) \\
  \bv &\mapsto \cP (\Ih^{p}\bv ).
\end{align*}

\begin{figure}
 \begin{center}
  \begin{subfigure}[b]{0.3\textwidth}
    \begin{center}
      \includegraphics[width=\textwidth]{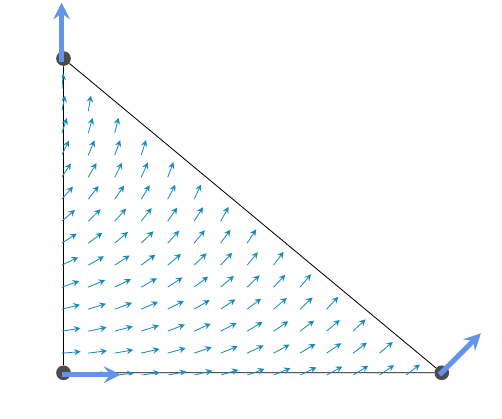}
    \end{center}
      \caption{Smooth}
  \end{subfigure} 
   \begin{subfigure}[b]{0.3\textwidth}
    \begin{center}
      \includegraphics[width=\textwidth]{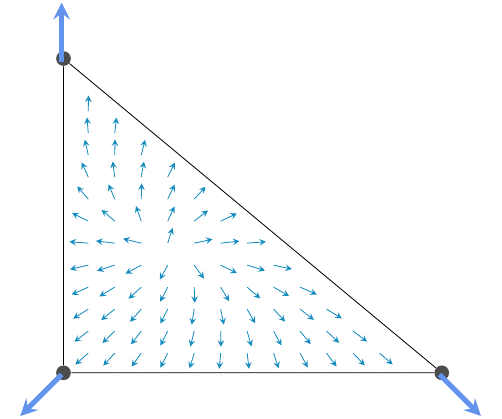}
    \end{center}
       \caption{Point singularity}
   \end{subfigure} 
   \begin{subfigure}[b]{0.3\textwidth}
    \begin{center}
      \includegraphics[width=\textwidth]{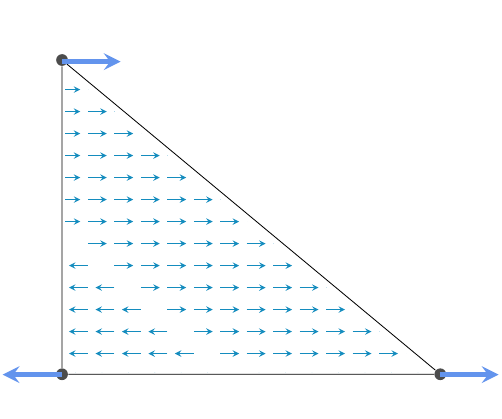}
    \end{center}
       \caption{Line singularity}
   \end{subfigure}
  \end{center}
  \caption{Three first-order projection-based finite element functions on a triangle,
  with values in $S^1$.
  The prescribed values at the Lagrange points are enlarged for better visibility.}
  \label{fig:gfe_functions}
 \end{figure}

Similar to the nonconforming discretization we then introduce the subspace of functions
that comply with the boundary condition at the boundary Lagrange points
\begin{align*}
  \cA_{h}^\textup{proj}
  \colonequals
  \Big\{ \bv_h  \in \cS^{p,\textup{proj}}(\Omega;S^{m-1}) \; :\;
   \bv_h(z) = \bu_\DD(z) \text{ for all } z \in \mathcal{L}_h \cap \partial \Omega
  \Big\}.
\end{align*}
Note that the correspondence between functions in $\cS^{p,\textup{proj}}$ and their sets
of values $c \in (S^{m-1})^N$ (with $N \colonequals \abs{\mathcal L_h}$) at the Lagrange points
is one-to-one. This allows to treat projection-based finite elements algorithmically as such sets
of values, just as in the case of standard finite elements.

It was shown in \cite{GFEprojection} that elements of $\cS^{p,\textup{proj}}$
are in $H^1(\Omega; S^{m-1})$
as long as configurations are avoided where the projection $P(\xi) \colonequals \xi / \abs{\xi}$
becomes undefined.  By this embedding property, the subset of continuous
projection-based finite elements (with fixed boundary conditions) decomposes into
disconnected homotopy classes just like $C(\Omega;S^{m-1})$ does.

However, as a particularity, the space $\cS^{p,\textup{proj}}(\Th;S^{m-1})$ also
contains maps with singularities. These appear whenever the values at the
Lagrange nodes on an element $T$ are such that $T$ contains a zero of the
Lagrange interpolation in $\R^m$, because there, $P : \R^m \to S^{m-1}$
is not defined. Figure~\ref{fig:gfe_functions} shows two such configurations
for the case $m=n=2$. It is easily seen that finite element maps with singularities can be classified
according to the dimension of the support of the singularity.  For example,
on a triangle, the Lagrange interpolation can be zero in a single point or on
a line intersecting the triangle.  Consequently, the projection-based finite element
functions can have point or line singularities.  However, as the values at the Lagrange points
are constrained to have unit length, the singularities cannot appear everywhere within $T$.
In particular, singularities at Lagrange points are not possible.

The precise regularity of these singular finite element functions remains
an open question.
The line-singular function on the right of Figure~\ref{fig:gfe_functions}
is discontinuous on a line, and therefore clearly not an element of $W^{1,q}$
(because it is not absolutely continuous on almost every line parallel to the coordinate directions~\cite{evans2018measure}).
However, the precise Sobolev regularity of point-singular finite element functions has not
been worked out yet.

Depending on the viewpoint, the singular projection-based finite element functions
are a nuisance or a feature. If the aim is to approximate smooth functions,
then the appearance of singular finite element functions can be problematic,
because one is leaving the realm of established error theory, and implementations
have to guard against division by zero.
Luckily, this happens only rarely in practice, usually when the grid is too coarse.

When trying to approximate singular functions, on the other hand, having singular approximation functions
at hand may be considered an asset. Indeed, singular functions are
difficult to approximate in standard finite element spaces, and introducing
additional bespoke singular finite element functions (as, e.g., in XFEM~\cite{fries_belytschko:2010})
is a standard way to increase the approximation power. As we will see below,
for harmonic maps we do not observe a clear advantage, and the topic needs
further investigation.
An additional technical problem is that special quadrature rules are required
to integrate the singular functions. We have experimented with Gauß--Legendre
und low-order composite rules, and did not find one to be clearly superior
to the other. Ideally, some sort of adaptive rule would be used.

\medskip  

As for the nonconforming discretization,
test functions are constructed as variations of finite element functions.
As functions in $\cS^{p,\textup{proj}}$ map into $S^{m-1}$ everywhere, such variations are
vector fields $\bvarphi_h : \Omega \to \R^m$
for which $\bu_h(x) \cdot \bvarphi_h(x) = 0$ for all $x\in\Omega$.%
\footnote{Here we tacitly interpret tangent vectors of $S^{m-1}$ as elements of the
surrounding space $\R^m$.}
Using a construction similar to the one given in \cite{sander2016test}, one can show that
the test function space at a function $\bv_h\in \cS^{p,\textup{proj}}$ is isomorphic to
$\prod\limits_{z\in\cL_h} T_{\bv_h(z)}S^{m-1}$. To construct this isomorphism,
let $I^\textup{proj}$ be the operator that maps coefficient sets in $(S^{m-1})^N$
to functions in $\cS^{p,\textup{proj}}$ and let $\bv_h\in \cS^{p,\textup{proj}}$ with coefficients
$v_1,\dots, v_N$.
Then, for all $b_i \in T_{v_i}S^{m-1}$, $i=1,\dots, N$ and $x \in \Omega$,
the corresponding test function~$\bvarphi_h$ can be evaluated by
 \begin{align*}
  \bvarphi_h[b_1,\dots,b_N](x)
       & = \sum_{i=1}^N \frac{\partial I^\textup{proj}(v_1,\dots,v_N;x)}{\partial v_i} \cdot b_i.
 \end{align*}
In line with Section~\ref{sec:discret-nonconf} we define $T_h^\textup{proj}(\bv_h)$
as the set of test functions at $\bv_h \in \mathcal{A}_h^\text{proj}$ that vanish on the boundary nodes.

As the admissible set $\cA_{h}^\textup{proj}$ is not a vector space,
the standard approximation error theory of finite elements does not apply.
However, for problems with smooth solutions there are rigorous optimal a priori
discretization error estimates
both for the $L^2$ norm and the $H^1$ norm \cite{GFE, GFEprojection}.
The proofs generalize results like the Bramble--Hilbert lemma
in order to prove optimal interpolation error estimates.
Then, the $H^1$ discretization error is estimated with a nonlinear version of the
Céa Lemma~\cite{grohs2015optimal}. To get bounds on the $L^2$ error,
a generalized Aubin--Nitsche lemma for predominantly quadratic energies such as the Dirichlet energy  is proved~\cite{harderingL2,hardering2017aubin,GFE}.
Here as well, the theoretical results are not restricted to the unit sphere, but hold for more general Riemannian manifolds~\cite{GFE}.

\section{Solvers for the discrete problems} \label{sec:solvers}

We now discuss two solver algorithms for the algebraic problems that arise from the discretizations
in Section~\ref{sec:discretizations}.
Both solvers presented here are based on energy minimization, and they are therefore unlikely
to find stationary harmonic maps that are not locally minimizing.
One of them is a conforming algorithm, i.e., each iterate is guaranteed to be an element
of the discrete admissible set. For the other one, this only holds in the limit
of vanishing step size.
We use the symbol $\cA_h$ to denote either the admissible set $\cA_h^\textup{nc}$
or $\cA_h^{\textup{proj}}$, and $T_h(\bv_h)$ to denote the corresponding test function space
at $\bv_h\in \cA_h$, when the difference is irrelevant.

\subsection{Nonconforming discrete gradient flow}
\label{sec:gradient-flow}

\begin{figure}
  \begin{center}
     \begin{tikzpicture}[scale=2, fatnode/.style={circle, fill=blue!50, inner sep=1.75pt}, tnode/.style={circle,fill=gray!50, inner sep=1.75pt} ]
          \coordinate (center) at (0,0);
      
          \draw[line width = 0.5mm] (center) circle (1.0);
      
          \node[fatnode] (A) at (0.71,-0.71) {};
          \node[fatnode] (B) at (1.35,0.0) {};
          \node[fatnode] (C) at (1.35,0.95) {};
          \node (F) at (0.6,1.70) {\textbf{\ldots}};
      
          \path[-Stealth] (center) edge node[left, pos=1.0, inner sep=6pt] {$\mathbf{u}_h^0$} (A);
          \path[-Stealth,color=gray] (A) edge node[right, pos=0.35]  {$\tau d_t\mathbf{u}_h^1$} (B);
          \path[-Stealth] (center) edge node[above,pos=0.87, inner sep=1pt ]  {$\mathbf{u}_h^1$} (B);
          \path[-Stealth,color=gray] (B) edge node[right, pos=0.5 ]  {$\tau d_t\mathbf{u}_h^2$} (C);
          \path[-Stealth] (center) edge node[below, pos=0.90, inner sep=1.5pt]  {$\mathbf{u}_h^2$} (C);
          \path[-Stealth,color=gray] (C) edge node[right, pos=0.6 ] {$\tau d_t\mathbf{u}_h^3$} (F);
      \end{tikzpicture}
 \caption{Steps of the discrete gradient flow}
 \label{fig:gradient_flow_steps}
\end{center}
\end{figure}

The first solver interprets the
minimization problem for the energy~\eqref{eq:dirichlet-energy} over $\cA_h^\textup{nc}$
as being posed on $(\R^m)^N$, $N \colonequals \abs{\mathcal{L}_h}$, subject to the constraint
$\abs{\bu(z)} = 1$ for each Lagrange point $z \in \mathcal{L}_h$.
Its approach to deal with the constraint is to linearize it and then,
starting from some initial configuration $\bu_h^0 \in \cA_h$, follow a discrete gradient flow for the energy~\eqref{eq:dirichlet-energy}, imposing the linearized constraint in the Lagrange points of the finite element space.
In order to define the discrete gradient flow, we use the $L^2$ scalar product of the gradients $(\,\nabla \cdot\,,\nabla \cdot\,)$.
Writing $P_{\bv_h}$ for the $H^1$-projection into the linear space $T_h(\bv_h)$
at a $\bv_h \in \mathcal{S}^p(\Th;\R^m)$, and introducing a pseudo-time variable $t$,
the gradient flow is
\begin{equation*}
 \frac{\partial \bu_h}{\partial t}
 =
 -P_{\bu_h}(\nabla_{H^1} E[\bu_h]),
\end{equation*}
with weak form
\begin{equation*}
 \Big( \nabla \frac{\partial \bu_h}{\partial t}, \nabla \bvarphi_h\Big)
 =
 -(\nabla \bu_h, \nabla \bvarphi_h).
\end{equation*}
The proposed algorithm is then an implicit Euler method for this flow.

\medskip

\begin{algorithm}[H]\label{alg:gradflow}
 \caption{Discrete gradient flow}
 \SetAlgoLined
 \DontPrintSemicolon
 \KwIn{initial iterate $\bu_h^0 \in \cA_h^\textup{nc}$, step size $\tau > 0$, stopping threshold $\epsilon_\mathrm{stop} > 0$}
 \For{$k=0,1,2,\dots$}
 {determine  $d_t \bu_h^{k+1} \in T_h(\bu_h^k)$ such that
  \[
    (\nabla d_t \bu_h^{k+1}, \nabla \bvarphi_h)
    =
    -(\nabla \bu_h^k, \nabla \bvarphi_h) - \tau (\nabla d_t \bu_h^{k+1}, \nabla \bvarphi_h)
    \text{ for all $\bvarphi_h \in T_h(\bu_h^k)$}
  \]
  \label{li:compute_tangent_update}

  set $\bu_h^{k+1} = \bu_h^k + \tau d_t\bu_h^{k+1}$ \\

  \If{$\abs{d_t \bu_h^{k+1}}_{H^1} \le \epsilon_\mathrm{stop}$}
  {
    \Return $\bu_h^{k+1}$
  }
 }
\end{algorithm}

\medskip

Determining the tangential correction in Line~\ref{li:compute_tangent_update} is a
linear elliptic equation on the vector space $T_h(\bu_h^k)$.
There are several ways how to enforce the restriction to the tangent space.
Following own previous work, our implementation uses Lagrange multipliers.
At each step~$k$, we solve the linear system
\begin{equation*}
  \begin{pmatrix}
    S & (A^k)^\top \\
    A^k & 0
  \end{pmatrix}
  \begin{pmatrix}
    d^{k+1} \\
    \lambda^{k+1}
  \end{pmatrix}
  =
  \begin{pmatrix}
    b \\
    0
  \end{pmatrix},
\end{equation*}
where $d^{k+1} \in (\R^m)^N$ is a representation of the discrete function $d_t \bu_h^{k+1}$
in the nodal basis of $\mathcal{S}^p(\mathcal{T}_h;\R^m)$, the matrix $S$
is the algebraic representation of the scalar product $(\nabla \cdot, \nabla \cdot)_{L^2}$
in the finite element space,
$b \in \R^{mN}$ encodes the explicit terms on the right hand side of Line~\ref{li:compute_tangent_update}
of Algorithm~\ref{alg:gradflow}, and $A^k$ is a matrix encoding the linearized
nodal constraints $d_t \bu_h^{k+1}(z) \cdot \bu_h^k(z) = 0$, $z \in \mathcal{L}_h$, with the corresponding
Lagrange multiplier $\lambda^{k+1} \in \R^N$.
The resulting linear systems are symmetric but indefinite, and are solved
with a direct solver if $\dim \Omega = 2$, and with a GMRES solver if $\dim \Omega = 3$.

Rigorous convergence analysis for this algorithm for the case of the first-order nonconforming
discretization of Chapter~\ref{sec:discret-nonconf} appears in~\cite{Bartels05,Bartels16}.
To obtain convergence, the time step size $\tau$ has to tend to zero with the mesh size $h$.
Monotonicity of the method implies that the norms of the corrections $d_tu_h^k$ decay like $\mathcal{O}\big((\tau k)^{-\frac{1}{2}}\big)$.
In this case the successive violations of the constraint, shown in~Figure~\ref{fig:gradient_flow_steps},
are controlled by the energy of the initial iterate and the step size,
\begin{equation}
 \label{eq:constraint-estimate}
  \max_{k=0,1,2,\dots} \int_{\Omega} \cI^1_h \big( \abs[\big]{|\bu_h^k|^2 -1}\big) \dv{x}
  \le
  c\tau E[\bu_h^0].
\end{equation}
In particular, the maximum violation is independent of the number of iterations, and the bound does not depend on structural properties of the underlying triangulation.
This observation holds for all target manifolds that are given as level sets of suitably regular functions.

In order to maintain the optimal convergence of the discretization, 
a violation of order $\mathcal{O}(h)$ in the constraint has to be guaranteed, 
which requires choosing $\tau = \mathcal{O}(h)$. 
This however reduces the speed of the convergence of the iteration which is $\mathcal{O}((\tau N)^{-\frac{1}{2}})$.

On the other hand, note that the linear convergence of the constraint violation
potentially spoils the approximation error when a higher-order finite element
space is employed. As a remedy a higher-order approximation in time may be
used to define the discrete gradient flow. Recently, it has been shown that
using a second-order backward differentiation formula leads to quadratic
constraint consistency~\cite{akrivis_bartels_palus:2024}.

An earlier version of Algorithm~\cite{Alouges97} projected each tangential correction
back onto the sphere at each Lagrange point. That way, no algebraic constraint violation
was accumulated, and the algorithm remained conforming. However, such projections may lead
to energy increase. Energy decrease can be guaranteed for certain types of triangulations,
but the corresponding restrictions limit the applicability of the algorithm
in three-dimensional situations \cite{Bartels05}. We do not consider this variant
of the algorithm any further in this text.

\subsection{Riemannian trust-region method}\label{sec:riemannianTR}

For a method that does not violate the algebraic sphere constraints, we turn to the field of
optimization on manifolds~\cite{absil2009optimization}.  Here, the admissible set
$\cA_h\in\{\cA_h^\textup{nc},\cA_h^\textup{proj}\}$ is viewed
as the product manifold $(S^{m-1})^N$.
The discrete problem for finding minimizing harmonic maps then has
the form of a minimization problem for the Dirichlet energy~\eqref{eq:dirichlet-energy}
on that manifold.

\medskip
\begin{algorithm}\label{alg:riemannianTR}
  \caption{Riemannian trust-region method}
  \SetAlgoLined
  \DontPrintSemicolon
  \KwIn{Initial iterate $\bu_h^0 \in \cA_h$,
  initial trust-region radius $\Delta_0 > 0 $,
  tolerances $\beta_1 > \beta_2 > 0$ 
  and stopping threshold $\epsilon_\mathrm{stop} > 0$
  }
  \For{$k=0,1,2,\dots$}
  {
   Solve \eqref{eq:TRsubproblem} for  $\bvarphi_h\in T_{h}(\bu_h^k)$ \tcp*{Compute tangential correction}
   \uIf{$\abs{ \bvarphi_h}_{H^1} < \epsilon_\mathrm{stop} $}
   {
   \Return $\bu_h^k$
   }
  \Else{

   Evaluate $\rho_k$ from \eqref{eq:qmeasure} \tcp*{Estimate model quality}

  \uIf{$\rho_k > \beta_1 $}
  {
    $\bu_h^{k+1} = \Exp_{\bu_h^k}(\bvarphi_h)$ and
    $\Delta_{k+1} = 2\Delta_k$  \tcp*{`Very successful' step}
  }
  \uElseIf{$\rho_k > \beta_2 $}
  {
    $\bu_h^{k+1} = \Exp_{\bu_h^k}(\bvarphi_h)$ and $\Delta_{k+1} = \Delta_{k}$
    \tcp*{`Successful' step}
  }
  \Else{
    $\bu_h^{k+1} = \bu_h^k$ and $\Delta_{k+1} = \frac{1}{2} \Delta_{k}$
    \tcp*{`Unsuccessful' step}
  }
  }
  }
\end{algorithm}
\medskip

In order to solve the minimization problem on $\cA_h$ we employ a Riemannian trust-region method.
Such a method generalizes standard trust-region methods to objective functionals
defined on a Riemannian manifold~\cite{absil2009optimization}.
In the spirit of Newton's method, the general idea
is to consider a quadratic model of the objective functional around the current iterate.
The algorithm then computes a correction $\bvarphi_h^k$
in the tangent space of the current iterate $\bu_h^k\in \cA_h$.
However, as the model is assumed to be accurate only in a neighborhood of $0$ on
$T_h(\bu_h^k)$, the tangential correction is restricted to remain in a ball
around~$0$ (the name-giving trust region), the radius of which is controlled adaptively.
To be more precise,
in each step $k$ the quadratic trust-region subproblem at an iterate $\bu_h^k$ is
\begin{align}
  \min_{\bvarphi_h\in T_h(\bu_h^k)} m_{k}(\bvarphi_h) &\colonequals
  E[\bu_h^k] + \left\langle \operatorname{Grad}
   E[\bu_h^k], \bvarphi_h \right\rangle
   +\frac{1}{2} \left\langle \operatorname{Hess}
    E[\bu_h^k](\bvarphi_h), \bvarphi_h \right\rangle \label{eq:TRsubproblem}
\end{align}
subject to
\begin{equation*}
  \norm{\bvarphi_h}_{\bu_h^k} \leq \Delta_k,
\end{equation*}
where $\norm{\cdot}_{\bu_h^k}$ is a suitable norm, and $\Delta_k$ is the current trust-region radius.
The terms $\operatorname{Grad}E$ and $\operatorname{Hess}E$ denote the Riemannian gradient and Hessian, respectively.
These can be computed from simple modifications of the corresponding Euclidean
quantities of the functional extended into the surrounding Euclidean space~\cite{absil2013extrinsic},
which take a particularly simple form for the unit sphere as a target manifold.
In a finite element context, one can conveniently solve the quadratic minimization
problems~\eqref{eq:TRsubproblem} by choosing the infinity norm $\norm{\cdot}_{\bv_h^k} = \norm{\cdot}_\infty$
with respect to some coordinate system to define the trust-region,
and then use a monotone multigrid (MMG) method for the minimization.
Detailed descriptions are given in~\cite{sander2016numerical,Sander12,kornhuber1997adaptive}.

Given an appropriate tangential update $\bvarphi_h^k$, the next iterate is then
$\bu_h^{k+1} = \exp_{\bu_h^k}(\bvarphi^k_h)\in \cA_h$, where $\exp_{\bu_h^k}$ is the
Riemannian exponential map of $\cA_h$ at $\bu_h^k$.
As $\cA_h$ is a product manifold, the exponential map
acts on each factor space $S^{m-1}$ separately. Conveniently, the exponential map of the sphere
has the closed-form expression
\begin{align*}
 \exp_\xi v
 =
 \cos \abs{v} \cdot \xi + \frac{\sin \abs{v}}{\abs{v}}\cdot v
 \qquad
 \text{where $\xi\in S^{m-1}$ and $v\in T_\xi S^{m-1}$}.
\end{align*}
Like the gradient-flow method of the previous section, the trust-region method computes
an update in a tangent space. However, while the gradient-flow algorithm uses Lagrange
multipliers to enforce tangentiality, the trust-region method
introduces a basis for $T_h(\bu_h^k)$ and solves the quadratic minimization problem
\eqref{eq:TRsubproblem} with respect to this basis.
As a consequence, the update problem involves fewer variables, and it is elliptic
instead of a saddle-point problem.
For coefficients in $S^2$, the basis is constructed as the push-forward of the canonical basis
of $\R^2$ under the inverse stereographic projection, which is conformal.
While the tangent problem of the gradient-flow method
could be easily formulated in terms of a basis for $T_h(\bu_h^k)$, the trust-region
problem involves an additional inequality constraint. Solving
minimization problems with inequality and equality constraints is possible,
but much more complicated.

Whether the new iterate $\bu_h^{k+1} = \Exp_{\bu_h^k}(\bvarphi_h)$ is actually accepted depends on whether it realizes sufficient energy decrease.
This is measured by the ratio of actual and expected energy decrease
\begin{align}
  \rho_k \colonequals
  \frac{E\big[\Exp_{\bu_h^k}(0)\big] - E\big[\Exp_{\bu_h^k}(\bvarphi_h)\big]}{m_{k}(0) - m_{k}(\bvarphi_h)}.
  \label{eq:qmeasure}
\end{align}
This quantity also controls the evolution of the trust-region radius~\cite{conn_gould_toint:2000}.

The global convergence of this method to first-order stationary points has been proved by Absil et al.~\cite{absil2009optimization}.
Moreover, depending on the accuracy of the inner solver, the local convergence is superlinear or even quadratic.
We refer to~\cite{absil2009optimization} for more details.
More applications of the method to manifold-valued finite element problems
appear in~\cite{sander2016numerical,nebel_sander_birsan_neff:2023}.

\section{Benchmarks problems}\label{sec:benchmarks-problems}

In this section we present the benchmark problems that will serve
to test the discretizations and solver algorithms.
All problems are posed on the open $n$-dimensional square $\Omega = (-\frac{1}{2}, \frac{1}{2})^n$.

\subsection{Inverse stereographic projection}

The first test problem constructs a harmonic map that is $C^\infty$ in the domain.
\begin{problembox}
 \begin{problem}[Inverse stereographic projection]\label{prob:inv_stereo}
   Let $\Omega = (-\frac{1}{2}, \frac{1}{2})^2$ and set $m=3$, i.e., the image space is $S^2\subset \R^3$.
   Denote the stereographic projection (from the north pole) by
   \begin{align*}
    \pi_{\rm st} &: S^2 \setminus \big\{(0,0,1)^T \big\} \to \R^2
    \qquad \qquad
    \pi_{\rm st}(\xi_1,\xi_2,\xi_3) \colonequals \left(\frac{\xi_1}{1-\xi_3},\frac{\xi_2}{1-\xi_3}\right)^T,
   \end{align*}
   and note that it is invertible with
   \begin{equation*}
     \pi_{\rm st}^{-1}(x_1,x_2)
     =
     \big(x_1^2 + x_2^2 + 1 \big )^{-1} \begin{pmatrix} 2 x_1 \\ 2x_2 \\ x_1^2 + x_2^2 - 1 \end{pmatrix}.
   \end{equation*}
   Find a harmonic map $\bu : \Omega \to S^2$ such that $\bu = \pi_{\rm st}^{-1}$ on the domain boundary $\partial\Omega$.
\end{problem}
\end{problembox}

Direct calculations show that $\pi_\text{st}^{-1}$ itself is a stationary map,
and it is in $C^\infty$. Its minimization properties
are discussed in~\cite{BreCor83}.

\subsection{Radial projection}
\label{sec:radial_projection}

The next two problems compute harmonic maps with a single singularity.
Both use the radial projection map
\begin{align*}
 \bu_\odot : \Omega \to S^{\dim \Omega-1},
 \qquad
 \bu_{\odot}(x) \colonequals \frac{x}{\abs{x}}
\end{align*}
as the boundary condition function.

\begin{problembox}
\begin{problem}[Radial projection]\label{prob:singular} \mbox{}
 \begin{enumerate}[label={(\alph*)},left=0pt]
  \item Find a harmonic map on $\Omega \colonequals \big(-\frac{1}{2},\frac{1}{2}\big)^3$
  with image in $S^2$ such that $\bu = \bu_\odot$ on $\partial \Omega$.

  \item Find a harmonic map on $\Omega \colonequals \big(-\frac{1}{2},\frac{1}{2}\big)^2$
  with image in $S^1$ such that $\bu = \bu_\odot$ on $\partial \Omega$.
 \end{enumerate}
\end{problem}
\end{problembox}

Both problems have $\bu_\odot$ as a solution.
A direct computations in polar coordinates shows that $\bu_\odot$ is an element
  of $W^{1,q}(\Omega;S^{m-1})$ for $q<n$ (see also \cite[Chap.\,3.2]{helein2008harmonic})
  and that it also belongs to $W^{2,q}(\Omega;S^{m-1})$ for $q<\frac{n}{2}$.
Further properties depend on the space dimensions:
\begin{enumerate}[label={(\alph*)}]
    \item \label{subprob:radial33}
      If $n = m = 3$, the radial projection $\bu_{\odot}$ is the
      unique minimizer of the Dirichlet energy in $H^1(\Omega,S^2)$. This was proven by Brezis, Coron and Lieb~\cite[Corollary~7.9]{BreCorLie86},
      based on the work of Schoen and Uhlenbeck~\cite{schoenUhlenbeck1982regularity}.
      The corresponding minimal energy on $\Omega$ is
      \[
      E[\bu_\odot] = 6 \int\limits_{-\pi/4}^{\pi/2} \frac{\pi - 2 \arctan(\frac{1}{\sin(\vartheta)})}{\sin(\vartheta)} \dv{\vartheta} \approx 7.674124,
      \]
      which can be verified with a computer algebra system.
      As $2 \nless \frac{n}{2}$ if $n=3$, we have $\bu_\odot \notin H^2(\Omega;S^2)$
      in this case.

    \item \label{subprob:radial22}
      If $n=m=2$, the radial projection $\bu_{\odot}$ is not even in $H^1(\Omega;S^1)$
      (see also ~\cite[p.\,71]{struwe2000variational}),
      and therefore $E[\bu_\odot]$ is not well-defined. However, it is still an element
      of $L^2(\Omega;S^1)$.
\end{enumerate}

The radial projection $\big(-\frac 12,\frac 12\big)^2 \to S^1$ is still harmonic
in the distributional sense, even though its Dirichlet energy is undefined.
We give a short proof of this result, which we have not found elsewhere.
\begin{lemma}
\label{lem:distributional_harmonic_map}
 If $n = m = 2$ the radial projection $\bu_{\odot} : x \mapsto \frac{x}{\abs{x}}$
 is harmonic in the distributional sense, i.e., $(\bu_\odot, \Delta \bvarphi) = 0$
 for any $\bvarphi \in C_0^\infty(\Omega;\R^m)$ with $\bu_\odot \cdot \bvarphi = 0$ a.e.
\end{lemma}

\begin{proof}
One directly checks that $\partial_i \bu_\odot$, $i=1,2$, and $-\Delta \bu_\odot$ are, respectively,
perpendicular and parallel to $\bu_\odot$, which itself coincides with the unit normal $\bn$
to the ball $B_\varepsilon(0)$ of radius $\varepsilon$ around the origin.
With this, a splitting of the integral and two integrations by parts lead to
\begin{multline*}
 \int_\Omega \bu_\odot \cdot \Delta \bvarphi \dv{x}
 =
  \int_{B_\varepsilon(0)} \bu_\odot \cdot \Delta \bvarphi \dv{x}
  + \int_{\partial(\Omega \setminus B_\varepsilon(0))} \bu_\odot \cdot (\nabla \bvarphi \bn) \dv{s} \\
  - \int_{\partial(\Omega \setminus B_\varepsilon(0))} (\nabla \bu_\odot \bn) \cdot \bvarphi \dv{s}
   + \int_{\Omega \setminus B_\varepsilon(0)} \Delta \bu_\odot \cdot \bvarphi \dv{x},
\end{multline*}
for any $\bvarphi\in C_0^{\infty}(\Omega;\R^m)$.
The first integral on the right-hand side is bounded in terms of $\varepsilon^2$,
the second integral is bounded in terms of $\varepsilon$,
the third integral vanishes since $\nabla \bu_\odot \bn = 0$,
and the fourth integral vanishes because $\bu_\odot$ is harmonic outside of $B_\varepsilon(0)$.
\end{proof}

\begin{figure}
  \centering
  \begin{subfigure}[b]{0.32\textwidth}
    \centering
    \includegraphics[width=\textwidth]{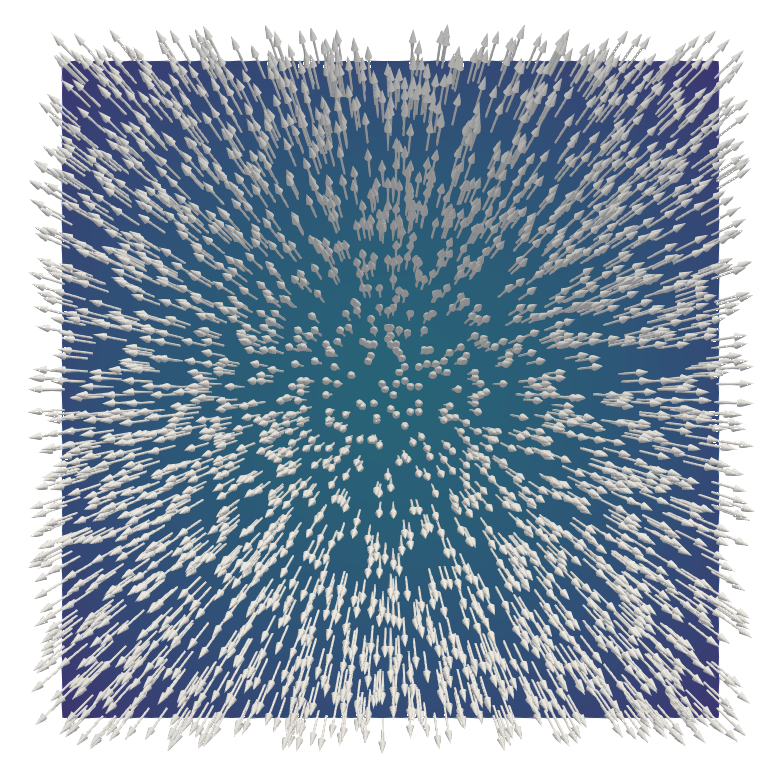}\\
    \caption{$\R^2 \supset \Omega \to S^2$ }
    \label{fig:radial-a}
  \end{subfigure}
  \begin{subfigure}[b]{0.32\textwidth}
    \centering
    \includegraphics[width=\textwidth]{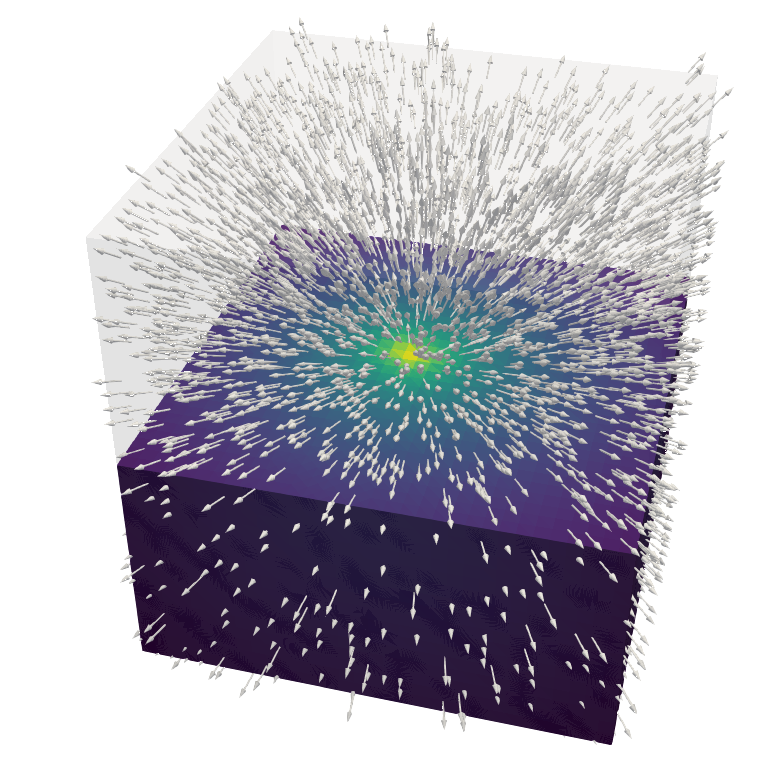}\\
    \caption{$\R^3 \supset \Omega \to S^2$ }
    \label{fig:radial-b}
  \end{subfigure}
  \begin{subfigure}[b]{0.32\textwidth}
    \centering
    \includegraphics[width=\textwidth]{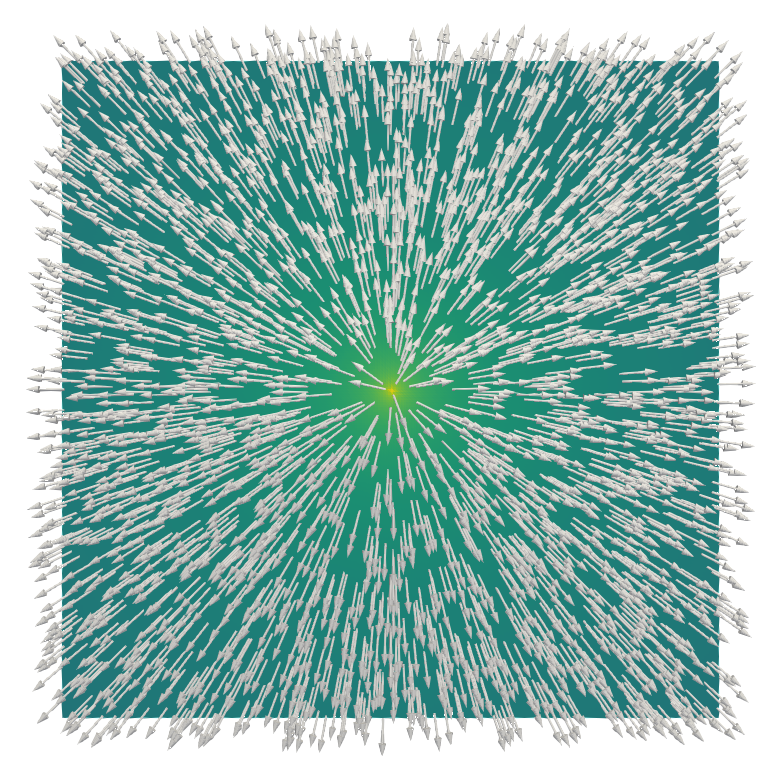}\\
    \caption{$\R^2 \supset \Omega \to S^1$ }
    \label{fig:radial-c}
  \end{subfigure}

  \caption{Discrete solutions $\bu_h$ for  Problems~\ref{prob:inv_stereo}, \ref{prob:singular}\ref{subprob:radial33} and \ref{prob:singular}\ref{subprob:radial22} }
  \label{fig:radial}
\end{figure}

\subsection{Multiple singularities}

\begin{figure}
  \centering
  \begin{subfigure}[b]{0.24\textwidth}
    \centering
    \includegraphics[width=\textwidth]{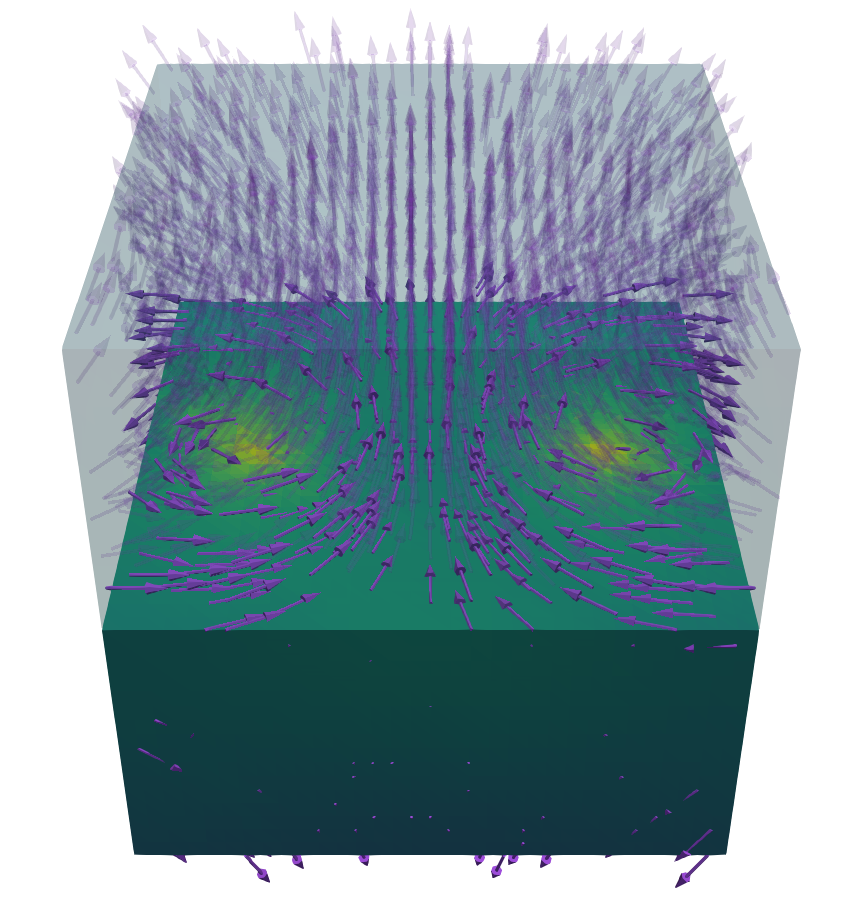}
  \end{subfigure}
  \begin{subfigure}[b]{0.24\textwidth}
    \centering
    \includegraphics[width=\textwidth]{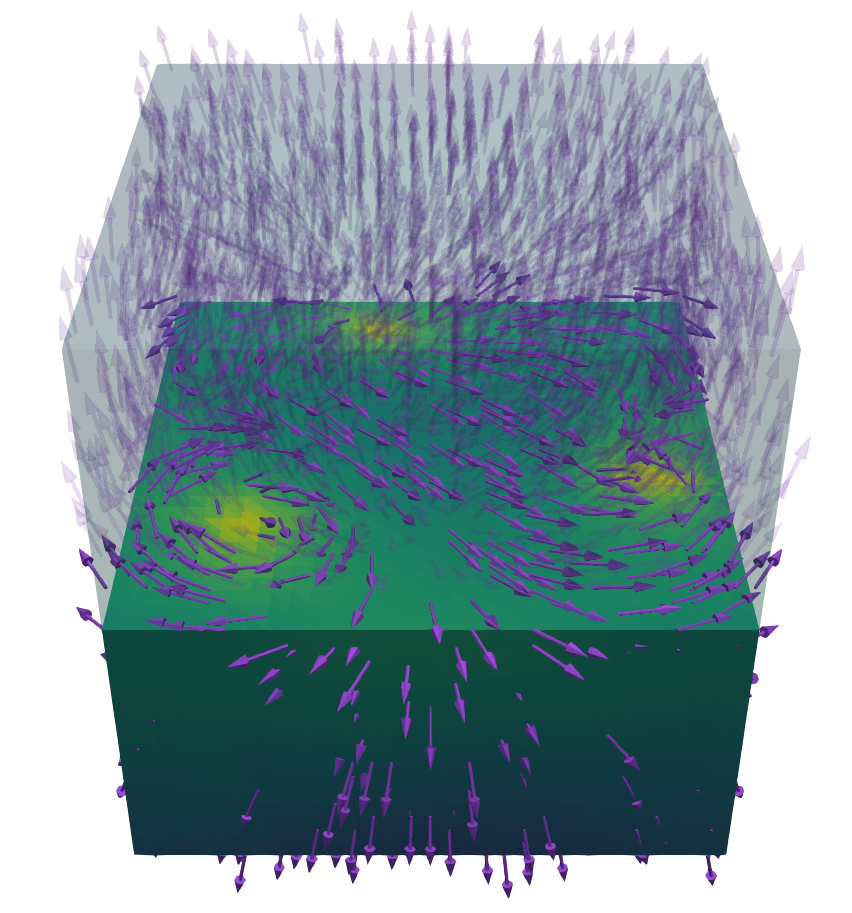}
  \end{subfigure}
  \begin{subfigure}[b]{0.24\textwidth}
    \centering
    \includegraphics[width=\textwidth]{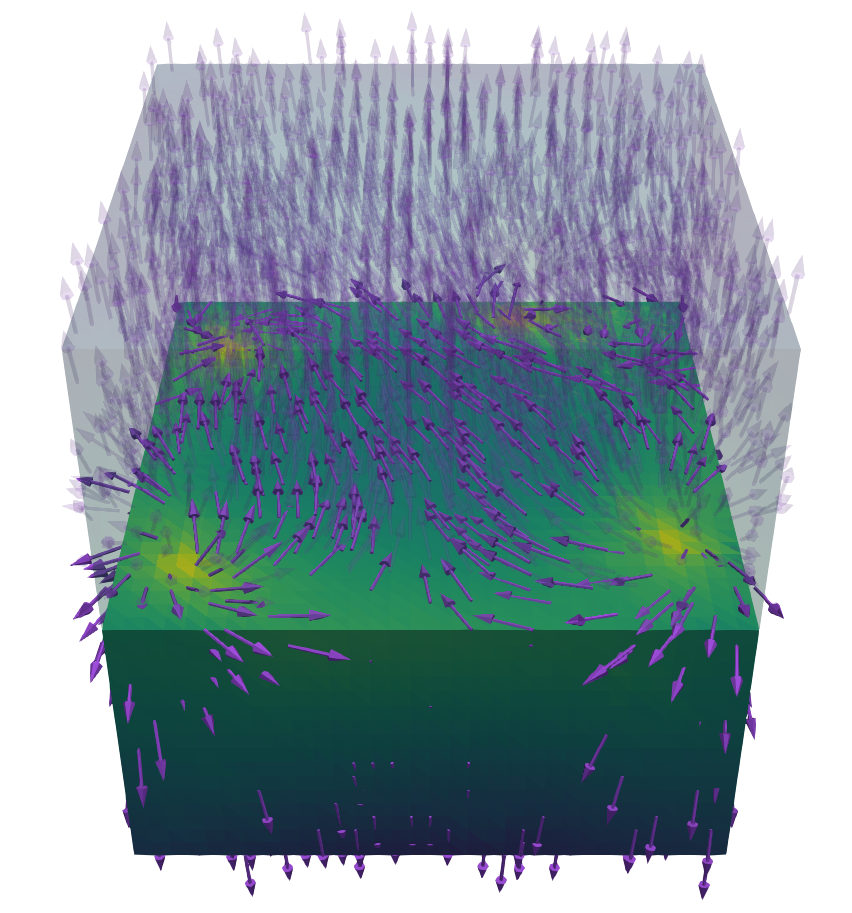}
  \end{subfigure}
  \begin{subfigure}[b]{0.24\textwidth}
    \centering
    \includegraphics[width=\textwidth]{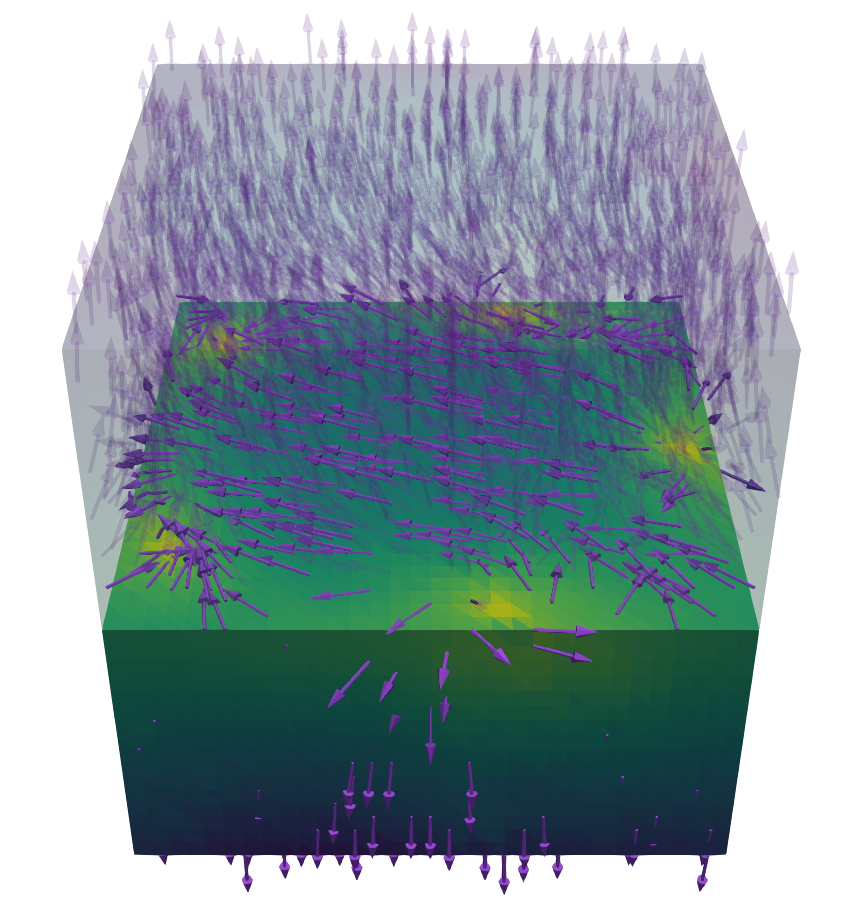}
  \end{subfigure}
  \caption{Stationary configurations of harmonic maps
  from $\big(-\frac{1}{2},\frac{1}{2}\big)^3$ to $S^2$ with multiple singularities,
  obtained with the nonconforming discretization and
  Algorithm~\ref{alg:gradflow} with initial singularity degrees $\kappa = 2,3,4,5$ (left to right).
  Singularities are located on the horizontal mid-surface, which is colored by the Frobenius norm of the discrete solution gradients.}
  \label{fig:highsing33}
\end{figure}

The final benchmark investigates harmonic maps with multiple singularities.
It also appears in~\cite{Alouges97}.

Recall that if a map $\bv \colon \Omega \to S^2$ is continuous in a ball $B_\varepsilon(x)$ around a point $x \in \Omega$ except at $x$ itself, the degree of the singularity at $x$ is defined as the topological winding number of $\bv|_{\partial B_\varepsilon(x)}$ with respect to~$x$, see~\cite{BreCorLie86}.

\begin{problembox}
\begin{problem}[Multiple singularities in $\R^3$]\label{prob:highsing}
  Interpret the stereographic projection $\pi_{\text{st}}$ as a map into the complex plane,
  and define
 \begin{equation}
 \label{eq:singularity_kappa}
  \uD^\kappa(x)
  \colonequals \pi_{\text{st}}^{-1} \circ g_\kappa \circ  \pi_{\text{st}} \circ \bu_\odot,
 \end{equation}
 where
 \begin{equation*}
   g_\kappa : \mathbb{C} \to \mathbb{C},
   \qquad
   g_\kappa(z)  \colonequals z^\kappa,
   \qquad
   \forall \kappa \in \mathbb{N}.
 \end{equation*}
 For given $\kappa \in \mathbb{N}$,
 find a harmonic map $\bu : \Omega = \big(-\frac 12,\frac 12\big)^3 \to S^2$
 such that $\bu = \bu_D^\kappa$ on $\partial \Omega$.
\end{problem}
\end{problembox}

The map $\bu_D^\kappa$ has a degree-$\kappa$ singularity at the origin.
For $\kappa =1$ this is identical to Problem~\ref{prob:singular}\ref{subprob:radial33}.
However, according to Brezis, Coron, and Lieb~\cite{BreCorLie86} (locally) minimizing harmonic maps
cannot have singularities of (absolute) degree greater than one.
Hence, the maps $\uD^\kappa$ are unstable configurations for $\kappa > 1$, and will
typically not be observed in simulations using energy-descent-based solver algorithms.
Figure~\ref{fig:highsing33} shows configurations for the cases $\kappa = 2,\dots,5$
obtained by numerical simulations using the nonconforming discretization and
the discrete gradient-flow solver. One can see that the initial
degree~$\kappa$ singularity of $\bu^\kappa_D$ splits up, and the system runs
into a stable state with $\kappa$ isolated singularities of unit degree each.

\section{Benchmarking the discretizations}
\label{sec:discretization_benchmarks}

We now numerically compare the two discretizations
of Section~\ref{sec:discretizations} on
the benchmark Problems~\ref{prob:inv_stereo}, \ref{prob:singular}\ref{subprob:radial33}
and~\ref{prob:singular}\ref{subprob:radial22} defined in Section~\ref{sec:benchmarks-problems}.
We omit Problem~\ref{prob:highsing},
which involves harmonic maps with multiple singularities.
The configurations of Problem~\ref{prob:highsing} are essentially multiple copies
of the single-singularity case of Problem~\ref{prob:singular}, and we
therefore do not expect to see any new effects.
Also, this problem has no closed-form solution, but solving it numerically
is also problematic (see Chapter~\ref{sec:solver_multiple_singularities} below).
As a consequence, we do not have a good reference to test the discretization errors with.

\subsection{Benchmarking procedure}

For the numerical experiments we employ uniform triangulations $\Th$ of $\Omega \subset \R^n$ into simplices.
Starting from triangulations with only 2 triangles (if $\Omega$ is two-dimensional) or
6 tetrahedra (if $\Omega$ is three-dimensional) we obtain sequences of test grids
by uniform refinement. We label these grids by their refinement level $r\in \mathbb{N}$.
A grid for $\Omega = \big(-\frac12,\frac12\big)^n$ with refinement level $r$ has
an element diameter of $h = 2^{-r}\sqrt{n}$.

We solve the algebraic systems resulting
from the different discretizations to machine precision
with the Riemannian trust-region solver of Section~\ref{sec:discret-conf}.
That way, no algebraic constraint violation occurs, and the algebraic error introduced by the iterative nature of the solver remains negligible.
We use the parameters $\Delta_0 = \frac{1}{2}$, $\delta_1 = 0.9$ and $\delta_2 = 10^{-2}$
for the outer solver. The inner multigrid solver is set to iterate
until the $H^1$ seminorm of the correction drops below $10^{-10}$.
The initial iterates will be given below.
The implementation is a hand-written \CC\  code based on the \textsc{Dune} libraries \cite{bastian2021dune,sander2020dune}. The \texttt{dune-alugrid}~\cite{Alugrid}
extension module is used as the grid data structure, and the \texttt{dune-gfe} module%
\footnote{\url{https://gitlab.mn.tu-dresden.de/osander/dune-gfe/}}
provides the implementation of the projection-based finite elements.

Let $\bu$ denote the known exact solution of a problem, and let $\bu_h$ be a finite element
approximation on a grid of maximal edge length $h$.
We can directly compute the errors in the $L^2$-norm $\norm{\bu_h - \bu}_{L^2}$
and the $H^1$-seminorm $\abs{\bu_h - \bu}_{H^1}$.
On our hierarchy of grids obtained by uniform refinement, we can estimate
the corresponding orders of convergence by
\begin{equation*}
\mathrm{EOC}^{L^2}_h \colonequals \log_2 \left(
  \frac{\| \bu_{2h} - \bu\|_{L^2_{}}}
  {\| \bu_{h} - \bu\|_{L^2}}
\right)
\qquad
\text{and}
\qquad
\mathrm{EOC}^{H^1}_h \colonequals \log_2 \left(
  \frac{| \bu_{2h} - \bu |_{H^1_{}}}
  {| \bu_{h} - \bu |_{H^1}}
\right).
\end{equation*}

Note that for the nonconforming discretization, the approximate solutions $\bu_h$ are piecewise polynomial functions,
and therefore the harmonic energy can be integrated exactly. However, the reference
solutions $\bu$ are not polynomial, and therefore computing the errors and convergence
orders always involves a quadrature error.
The projection-based finite elements, however, are not piecewise polynomials themselves,
and no exact quadrature formula is known even for the harmonic energy.
Here, we
use Gauss--Legendre quadrature of second and sixth order for the harmonic energy and the EOCs, respectively.
Finding more appropriate quadrature formulas for geometric finite elements remains
an interesting research subject of its own.

\subsection{Inverse stereographic projection}\label{sec:radial23}
We begin with Problem~\ref{prob:inv_stereo}, i.e., the computation of a harmonic map
from $\Omega = \big(-\frac{1}{2}, \frac{1}{2}\big)^2$ to $S^2$, with boundary data
given by the restriction of the inverse stereographic projection $\pi_\text{st}^{-1}$ to $\Omega$.
One solution to this is the inverse stereographic projection itself,
which is smooth.

We measure the discretization errors for finite elements of orders $p=1$ and $p=2$,
starting the solver from the interpolant of the inverse stereographic projection.
The results are shown in Table~\ref{tab:inv_stereo}.
The rates agree with what would be expected for a linear problem, i.e.,
they are iin $\cO(h^{p+1})$ for the $L^2$-error and in $\cO(h^{p})$ for the $H^1$-error.
For the projection-based finite element discretization this has been proven in~\cite{GFEprojection,GFE}.
For nonconforming discretizations with $p=1$, optimal error estimates in the
energy norm have been derived for a corresponding saddle-point formulation in~\cite{huWinther2009saddle} for the case $m=n=2$,
and for more general two- and three-dimensional settings in~\cite{BaPaWa2023}.

Table~\ref{tab:inv_stereo} also shows the discrete energies of the minimizers for the different grids.
The values produced by the two different discretizations
roughly agree---more so for the second-order finite elements.
Observe how the minimal energy increases with increasing mesh refinement for the nonconforming discretization, whereas it decreases for the conforming discretization. The latter would be the expected behavior for nested finite element spaces.
However, neither the conforming nor the nonconforming finite elements form nested approximation space hierarchies with respect to the refinement $r$.

\begin{table}\footnotesize
  \addtolength{\tabcolsep}{-1.5pt}
  \footnotesize
  \begin{tabular}{c r | c c c | c c c}
    \multicolumn{2}{c}{} &\multicolumn{3}{c}{nonconforming} &\multicolumn{3}{c}{conforming}\\
    $r$ &$\abs{\Th}$
    &$E[\bu_h]$ &$\mathrm{EOC}^{L^2}_h$ & $\mathrm{EOC}^{H^1}_h$
    &$E[\bu_h]$ &$\mathrm{EOC}^{L^2}_h$ & $\mathrm{EOC}^{H^1}_h$ \\ \hline\hline
    & & \multicolumn{6}{c}{order $p=1$} \\
    1 &              8 & \num{2.66667}     &  -          &  -          & \num{3.19238}     &  -          &  -          \\
    2 &             32 & \num{2.91956}     & \ordernum{1.92991}   & \ordernum{0.94893}   & \num{3.06189}     & \ordernum{2.15187}   & \ordernum{1.04435}   \\
    3 &            128 & \num{2.98648}     & \ordernum{1.98377}   & \ordernum{0.98789}   & \num{3.02271}     & \ordernum{2.00503}   & \ordernum{1.00177}   \\
    4 &            512 & \num{3.00343}     & \ordernum{1.99598}   & \ordernum{0.99701}   & \num{3.01253}     & \ordernum{1.99653}   & \ordernum{0.99991}   \\
    5 &           2048 & \num{3.00768}     & \ordernum{1.99899}   & \ordernum{0.99925}   & \num{3.00996}     & \ordernum{1.99875}   & \ordernum{0.99994}   \\
    6 &           8192 & \num{3.00874}     & \ordernum{1.99975}   & \ordernum{0.99981}   & \num{3.00931}     & \ordernum{1.99967}   & \ordernum{0.99999}   \\
    7 &          32\,768 & \num{3.00901}     & \ordernum{1.99994}   & \ordernum{0.99996}   & \num{3.00915}     & \ordernum{1.99992}   & \ordernum{0.99999}   \\
    8 &         131\,072 & \num{3.00908}     & \ordernum{1.99999}   & \ordernum{0.99998}   & \num{3.00911}     & \ordernum{1.99998}   & \ordernum{1.000}       \\
    \hline
    & & \multicolumn{6}{c}{order $p=2$} \\
    1 &              8 & \num{3.00386}     &  -              &  -              & \num{3.0295}      &  -          &  -         \\
    2 &             32 & \num{3.00877}     & \ordernum{2.96204}   & \ordernum{1.93833}   & \num{3.01052}     & \ordernum{2.92379}   & \ordernum{1.92721}  \\
    3 &            128 & \num{3.00908}     & \ordernum{2.98974}   & \ordernum{1.98676}   & \num{3.00919}     & \ordernum{3.02682}   & \ordernum{1.98442}  \\
    4 &            512 & \num{3.0091}      & \ordernum{2.99776}   & \ordernum{1.99702}   & \num{3.0091}      & \ordernum{3.01314}   & \ordernum{1.99627}  \\
    5 &           2048 & \num{3.0091}      & \ordernum{2.99948}   & \ordernum{1.99929}   & \num{3.0091}      & \ordernum{3.00382}   & \ordernum{1.99909}  \\
    6 &           8192 & \num{3.0091}      & \ordernum{2.99988}   & \ordernum{1.99983}   & \num{3.0091}      & \ordernum{3.00099}   & \ordernum{1.99977}  \\
    7 &          32\,768 & \num{3.0091}      & \ordernum{2.99997}   & \ordernum{1.99996}   & \num{3.0091}      & \ordernum{3.00025}   & \ordernum{1.99995}  \\
    \hline
  \end{tabular}%
  \caption{Experimental discretization error convergence orders for Problem~\ref{prob:inv_stereo}}
  \label{tab:inv_stereo}
\end{table}

\subsection{Radial projection}

We repeat the same experiment for Problems~\ref{prob:singular}\ref{subprob:radial33} and \ref{prob:singular}\ref{subprob:radial22},
which both lead to harmonic maps with a singularity.

\subsubsection{The case $n = m = 3$}\label{ex:radial33}

We first test the discretizations on Problem~\ref{prob:singular}\ref{subprob:radial33},
the solution of which is the radial projection $\bu_\odot : x \mapsto \frac{x}{\abs{x}}$ on the
three-dimensional domain $\Omega = \big(-\frac12,\frac12\big)^3$.
Recall that the radial projection map $\bu_\odot$ is in $H^1(\Omega;S^2)$,
but not in $H^2(\Omega;S^2)$ (more details in Chapter~\ref{sec:radial_projection}).

We cannot pick the solution as the initial iterate, and instead we choose
\begin{equation}
\label{eq:radial_projection_with_origin_tweak_3d}
 \bu_h^0 \in \cA_h,
 \qquad
 \bu_h^0(x)
 =
 \begin{cases}
  \bu_\odot(x) & \text{if $x \in \mathcal{L}_h$ and $x \neq (0,0,0)^T$}, \\
  (0,0,1)^T    & \text{if $x \in \mathcal{L}_h$ and $x = (0,0,0)^T$}.
 \end{cases}
\end{equation}
The special treatment for the value at $(0,0,0)^T \in \Omega$ is necessary
because it is a Lagrange point, but $\bu_\odot$ is not defined there.
(Even though projection-based finite elements can represent singular functions,
they cannot represent functions that are singular at Lagrange points.)
The choice of $(0,0,1)^T$ for the value at $x = (0,0,0)^T$ is arbitrary;
note that it breaks some of the problem's inherent symmetry.

The experimental results are listed in Table~\ref{tab:radial33}.
The measured orders of convergence are much lower now that the solution
is not in $H^2(\Omega;S^2)$, and they have a larger variance.
For both discretizations we obtain orders that are a little under $1$
for the $L^2$-error and around $0.4$ for the $H^1$-error.
Given the regularity of the solution $\bu_\odot$, these convergence orders are plausible.
However, no rigorous a priori discretization error bounds exist currently
for either discretization.  The conforming discretization seems to perform a bit better
than the nonconforming one, but this may be coincidental.
Its $L^2$-error is much closer to $1$,
but it has an outlier at $\approx 0.5$ for the grid with $r=3$.
The reason for this is unclear.

Again, Table~\ref{tab:radial33} also shows the energies of the minimizers for the two discretizations.
The energies increase as the grid is refined for the nonconforming discretization
(with one exception), whereas it decreases for the nonconforming discretization.
This is the same behavior as in the smooth case.

\begin{table}\footnotesize
  \addtolength{\tabcolsep}{-1.5pt}
  \footnotesize
  \begin{tabular}{c r | c c c | c c c}
  \multicolumn{2}{c}{}  &\multicolumn{3}{c}{nonconforming} &\multicolumn{3}{c}{conforming}\\
    $r$ &$\abs{\Th}$
    &$E[\bu_h]$ &$\mathrm{EOC}^{L^2}_h$ & $\mathrm{EOC}^{H^1}_h$
    &$E[\bu_h]$ &$\mathrm{EOC}^{L^2}_h$ & $\mathrm{EOC}^{H^1}_h$ \\ \hline\hline
    & & \multicolumn{6}{c}{order $p=1$} \\
    1&           48 & \num{5.30787}     & -                    & -                    & \num{7.85756}     & -        & -         \\
    2&          384 & \num{6.6103}      & \ordernum{1.05236}   & \ordernum{0.31193}   & \num{8.00068}     & \ordernum{1.0227}    & \ordernum{0.43508}   \\
    3&         3072 & \num{7.16368}     & \ordernum{0.77591}   & \ordernum{0.36417}   & \num{7.93672}     & \ordernum{0.55501}   & \ordernum{0.47906}   \\
    4&      24\,576 & \num{7.42324}     & \ordernum{0.71998}   & \ordernum{0.37591}   & \num{7.81968}     & \ordernum{0.86305}   & \ordernum{0.40884}   \\
    5&     196\,608 & \num{7.54911}     & \ordernum{0.80537}   & \ordernum{0.40557}   & \num{7.75119}     & \ordernum{0.9283}    & \ordernum{0.44082}   \\
    6&  1\,572\,864 & \num{7.61155}     & \ordernum{0.88834}   & \ordernum{0.43838}   & \num{7.71377}     & \ordernum{0.96248}   & \ordernum{0.46362}  \\
    \hline
    & & \multicolumn{6}{c}{order $p=2$} \\
    1 &             48 & \num{6.64166}     & -                    & -                  &      \num{7.92037}     & -         & -         \\
    2 &            384 & \num{7.16368}     & \ordernum{0.86067}   & \ordernum{0.36303} &      \num{7.83252}     & \ordernum{1.25831}   & \ordernum{0.48029}   \\
    3 &           3072 & \num{7.40946}     & \ordernum{0.89374}   & \ordernum{0.4143}  &      \num{7.75332}     & \ordernum{1.20542}   & \ordernum{0.49549}   \\
    4 &        24\,576 & \num{7.53837}     & \ordernum{0.9254}    & \ordernum{0.43786} &      \num{7.71341}     & \ordernum{1.13831}   & \ordernum{0.49568}   \\
    5 &       196\,608 & \num{7.60528}     & \ordernum{0.95429}   & \ordernum{0.45886} &      \num{7.69365}     & \ordernum{1.08008}   & \ordernum{0.49657}   \\
    \hline
  \end{tabular}%
  \caption{Experimental discretization error convergence orders for Problem~\ref{prob:singular}\ref{subprob:radial33},
  with the singularity at $(0,0,0)^T$ and using the initial iterate defined in~\eqref{eq:radial_projection_with_origin_tweak_3d}}
  \label{tab:radial33}
  \end{table}

\subsubsection{Moving the singularity off the Lagrange point}

The situation of the previous section is somewhat special: The singularity
of the exact solution is right on a Lagrange point for all grid refinements,
and it therefore cannot be represented even with a projection-based finite element function.

To study whether this has harmful consequences, we redo the previous experiment with a singularity
that is slightly shifted.  More precisely, we keep the domain
$\Omega = \big(-\frac12,\frac12\big)^3$, but we set the Dirichlet boundary values such that
\begin{equation*}
 \bu_{\odot,\delta}
 \colonequals
 \frac{x-\delta}{\abs{x-\delta}},
 \qquad
 \delta \colonequals 2^{-6}\Big(\frac13, \frac13, \frac13\Big)^T
\end{equation*}
becomes a solution.
With this choice of $\delta$, the singularity is not on a Lagrange point
for all grids that we test with, and we can therefore use it as initial iterate
without modification.

Table~\ref{tab:radial33-movedSing} shows the results.
One can see that the new position of the singularity does have an influence on the convergence
behavior. For both types of discretizations, the measured orders are a little higher.
On the other hand, the orders are still highly unstable from level to level,
and it is not possible to single out a fixed value as the definite order.
There is even one outlier where the order becomes negative. Fully understanding the behavior
of the discretizations in such a nonregular scenario certainly requires further work.

  \begin{table}\footnotesize
    \addtolength{\tabcolsep}{-1.5pt}
    \footnotesize
  \begin{tabular}{c r | c c c | c c c}
  \multicolumn{2}{c}{}  &\multicolumn{3}{c}{nonconforming} &\multicolumn{3}{c}{conforming}\\
    $r$ &$\abs{\Th}$
    &$E[\bu_h]$ &$\mathrm{EOC}^{L^2}_h$ & $\mathrm{EOC}^{H^1}_h$
    &$E[\bu_h]$ &$\mathrm{EOC}^{L^2}_h$ & $\mathrm{EOC}^{H^1}_h$ \\ \hline\hline
    & & \multicolumn{6}{c}{order $p=1$} \\
    1&  48           & \num{5.2718}      & -                  & -                  & \num{7.81431}     & -         & -        \\
    2&  384          & \num{6.58986}     & \ordernum{1.07349} & \ordernum{0.32093} & \num{7.98562}     & \ordernum{1.04752}   & \ordernum{0.43623}  \\
    3&  3072       & \num{7.15066}     & \ordernum{0.82019} & \ordernum{0.38512} & \num{7.97772}     & \ordernum{0.34919}   & \ordernum{0.35397}  \\
    4& 24\,576       & \num{7.41552}     & \ordernum{0.81485} & \ordernum{0.39118} & \num{7.83361}     & \ordernum{0.96216}   & \ordernum{0.4033}   \\
    5& 196\,608      & \num{7.54486}     & \ordernum{1.0274}  & \ordernum{0.53183} & \num{7.74717}     & \ordernum{1.46202}   & \ordernum{0.77826}  \\
    6& 1\,572\,864   & \num{7.60933}     & \ordernum{1.47208} & \ordernum{0.40729} & \num{7.71176}     & \ordernum{1.69798}   & \ordernum{0.54909}  \\
    \hline 
    & & \multicolumn{6}{c}{order $p=2$} \\
    1&  48           & \num{6.63771}     & -                  & -                  & \num{7.91146}     & -        & -     \\
    2&  384          & \num{7.16236}     & \ordernum{0.84915} & \ordernum{0.44235} & \num{7.82815}     & \ordernum{1.36573}  & \ordernum{0.52394}\\
    3&  3072       & \num{7.40277}     & \ordernum{0.91849} & \ordernum{0.46594} & \num{7.75116}     & \ordernum{1.53011}  & \ordernum{0.59184}\\
    4& 24\,576       & \num{7.53318}     & \ordernum{1.1707}  & \ordernum{0.51209} & \num{7.71235}     & \ordernum{1.78979}  & \ordernum{0.54576}\\
    5& 196\,608      & \num{7.60219}     & \ordernum{1.78332} & \ordernum{0.80619} & \num{7.69314}     & \ordernum{-0.25463} & \ordernum{0.96739}\\
    \hline
  \end{tabular}%
  \caption{Experimental discretization error convergence orders for Problem~\ref{prob:singular}\ref{subprob:radial33},
  with singularity at $\delta = 2^{-6}(\frac13,\frac13,\frac13)^T$}
  \label{tab:radial33-movedSing}
  \end{table}

\subsubsection{The case $n=m=2$}\label{ex:radial22}
In the situation of Problem~\ref{prob:singular}\ref{subprob:radial22} the solution
is even less regular.
The radial projection $\bu_\odot : \big(-\frac 12,\frac 12\big)^2 \to S^1$
is still a harmonic map, but it is so only
in the distributional sense (Lemma~\ref{lem:distributional_harmonic_map}).  As a consequence,
the Dirichlet energy of $\bu_\odot$ is not a finite number.
The problem is still well-defined in finite element spaces when using numerical quadrature,
but we cannot hope for convergence
in the $H^1$ sense.

To compute finite element solutions we use the same approach as for Problem~\ref{prob:singular}\ref{subprob:radial33}
with $n=m=3$.  The singularity is on a Lagrange point again, and we set the
initial iterate to the vector $(1,0)^T$ there.
The results are shown in Table~\ref{tab:radial22}.
For both discretizations we still observe convergence in the $L^2$-sense to $\bu_\odot$.
The orders are even a bit less oscillatory than in the three-dimensional situation.
Indeed, the orders of the conforming discretization are slightly better than
the orders of the nonconforming one, and both
seem to be increasing slowly with increasing mesh resolution.
All orders stay below~1,
and it is unclear whether the orders would approach a limit on even finer grids.
In similar experiments in \cite{huWinther2009saddle}, Hu, Tai, and Winther observed the same linear convergence.

However, in the $H^1$-sense there is no convergence at all.
Table~\ref{tab:radial22} shows the errors, and one can observe that they converge to a fixed value
not equal to zero.

\begin{table}\footnotesize
  \addtolength{\tabcolsep}{-1.5pt}
  \footnotesize
  \begin{tabular}{c r | c c c | c c c}
  \multicolumn{2}{c}{ } &\multicolumn{3}{c}{nonconforming} &\multicolumn{3}{c}{conforming}\\
    $r$ &$\abs{\Th}$
    &$E[\bu_h]$ &$\mathrm{EOC}^{L^2}_h$ & ${| \bu_{h} - \bu |_{H^1}}$
    &$E[\bu_h]$ &$\mathrm{EOC}^{L^2}_h$ & ${| \bu_{h} - \bu |_{H^1}}$ \\ \hline\hline
    & & \multicolumn{6}{c}{order $p=1$} \\
    1&  8       &   \num{5.17157}     & -                   & \num{4.301e+00}    & \num{5.52824}     & -                  & \num{5.498e+00}   \\
    2&  32      &   \num{6.86567}     & \ordernum{0.6522}   & \num{4.373e+00}    & \num{6.79358}     & \ordernum{0.67645} & \num{6.695e+00}   \\
    3&  128     &   \num{8.83366}     & \ordernum{0.62064}  & \num{4.483e+00}    & \num{8.79102}     & \ordernum{0.82384} & \num{6.842e+00}   \\
    4&  512     &   \num{10.9553}     & \ordernum{0.76066}  & \num{4.529e+00}    & \num{10.9424}     & \ordernum{0.87769} & \num{6.853e+00}   \\
    5&     2048 &   \num{13.1187}     & \ordernum{0.83174}  & \num{4.544e+00}    & \num{13.1143}     & \ordernum{0.89724} & \num{6.855e+00}   \\
    6&     8192 &   \num{15.2927}     & \ordernum{0.86815}  & \num{4.548e+00}    & \num{15.2905}     & \ordernum{0.91013} & \num{6.856e+00}   \\
    7& 32\,768  &   \num{17.4693}     & \ordernum{0.89005}  & \num{4.549e+00}    & \num{17.4677}     & \ordernum{0.92004} & \num{6.856e+00}   \\
    8& 131\,072 &   \num{19.6467}     & \ordernum{0.90506}  & \num{4.550e+00}    & \num{19.6452}     & \ordernum{0.92799} & \num{6.856e+00}   \\
    \hline
    & & \multicolumn{6}{c}{order $p=2$} \\
    1&  8       &  \num{7.43037}     & -                    & \num{3.636e+00}   & \num{10.1924}     & -                    & \num{3.850e+00}  \\
    2&  32      &  \num{9.45182}     & \ordernum{0.643}     & \num{3.763e+00}   & \num{12.2994}     & \ordernum{0.93259}   & \num{3.883e+00}  \\
    3&  128     &  \num{11.5885}     & \ordernum{0.78686}   & \num{3.805e+00}   & \num{14.468}      & \ordernum{0.92072}   & \num{3.885e+00}  \\
    4&  512     &  \num{13.7554}     & \ordernum{0.84479}   & \num{3.818e+00}   & \num{16.6432}     & \ordernum{0.92314}   & \num{3.886e+00}  \\
    5& 2048     &  \num{15.9303}     & \ordernum{0.87559}   & \num{3.822e+00}   & \num{18.8202}     & \ordernum{0.92941}   & \num{3.887e+00}  \\
    6& 8192     &  \num{18.1072}     & \ordernum{0.89492}   & \num{3.823e+00}   & \num{20.9976}     & \ordernum{0.93546}   & \num{3.887e+00}  \\
    7& 32\,768  &  \num{20.2846}     & \ordernum{0.90858}   & \num{3.823e+00}   & \num{23.1752}     & \ordernum{0.94072}   & \num{3.887e+00}  \\
    \hline
  \end{tabular}%
  \caption{Experimental discretization error convergence orders for Problem~\ref{prob:singular}\ref{subprob:radial22}}
  \label{tab:radial22}
  \end{table}

\section{Benchmarking the solvers}\label{sec:solver-benchmark}

\begin{figure}
  \centering
  \begin{subfigure}[b]{0.32\textwidth}
    \centering
    \includegraphics[width=\textwidth]{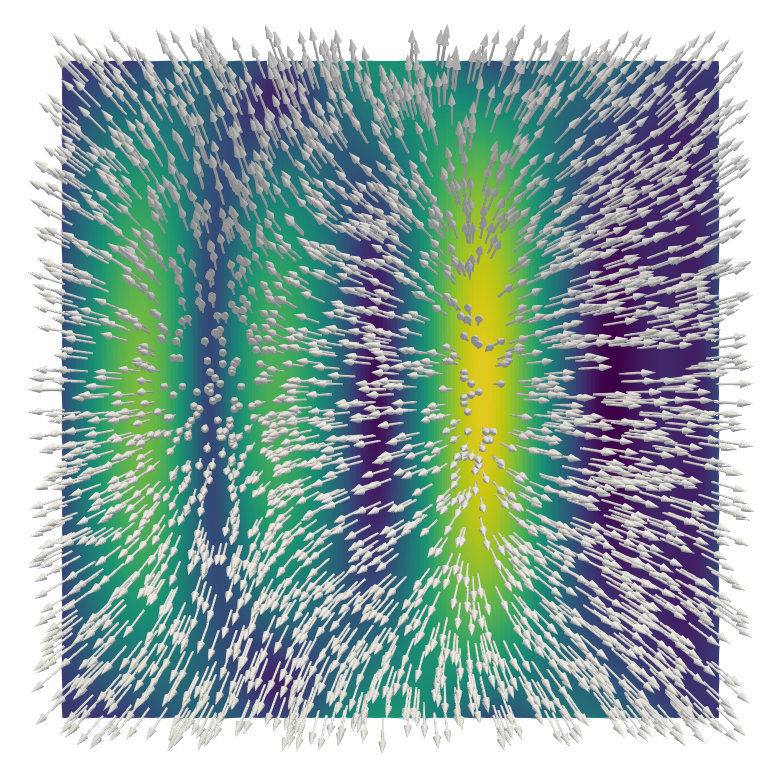}
    \caption{$\R^2 \supset \Omega \to S^2$ }
    \label{fig:radial-a-perturbed}
  \end{subfigure}
  \begin{subfigure}[b]{0.32\textwidth}
    \centering
    \includegraphics[width=\textwidth]{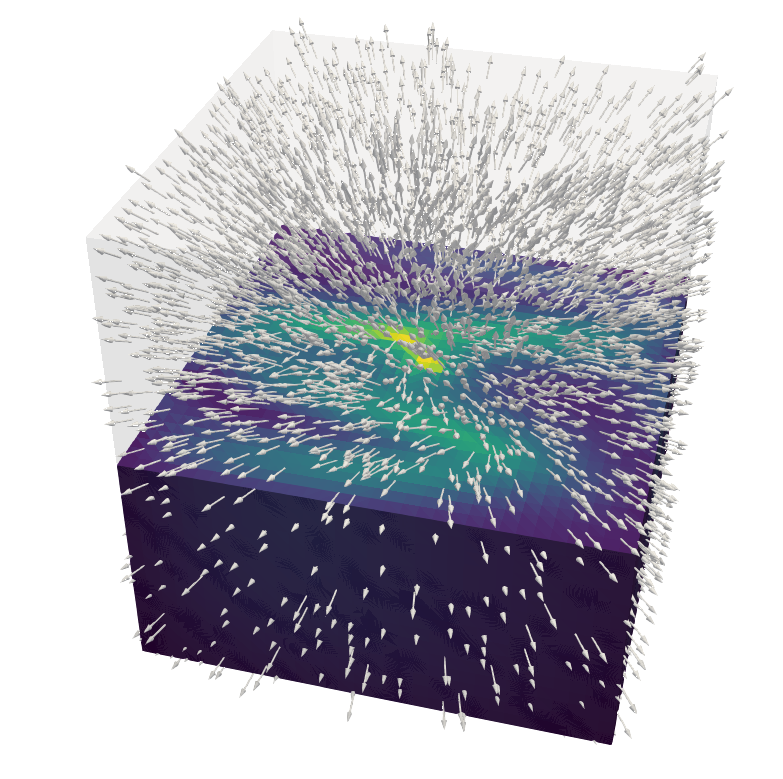}
    \caption{$\R^3 \supset \Omega \to S^2$ }
    \label{fig:radial-b-perturbed}
  \end{subfigure}
  \begin{subfigure}[b]{0.32\textwidth}
    \centering
    \includegraphics[width=\textwidth]{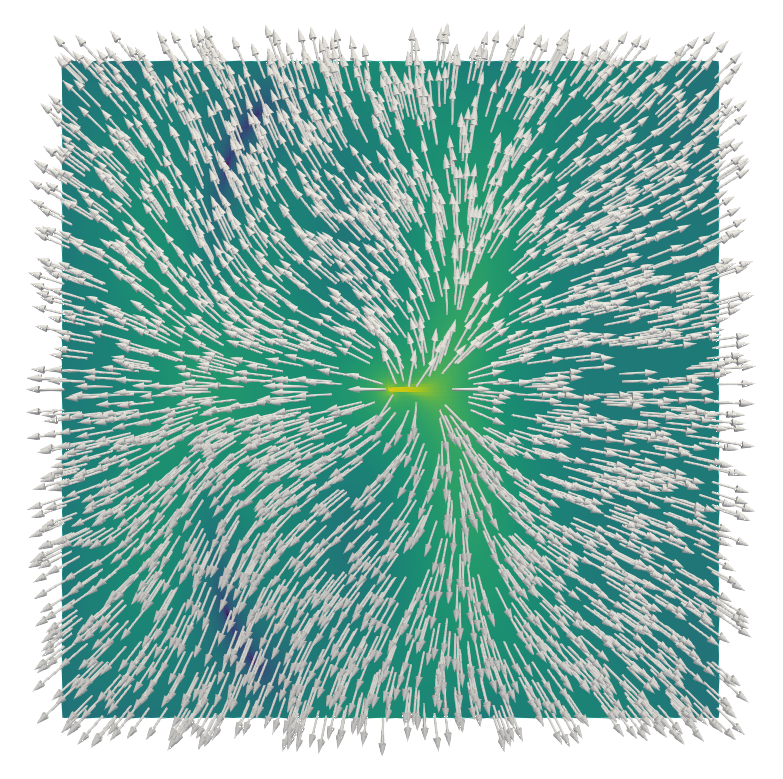}
    \caption{$\R^2 \supset \Omega \to S^1$ }
    \label{fig:radial-c-perturbed}
  \end{subfigure}

  \caption{Initial configurations $\bu_h^0$ for  Problems~\ref{prob:inv_stereo}, \ref{prob:singular}\ref{subprob:radial33}
   and \ref{prob:singular}\ref{subprob:radial22}. These are perturbations of the problem solutions. }
  \label{fig:radial-perturbed}
\end{figure}

In this second chapter of numerical tests we benchmark the solver algorithms of Chapter~\ref{sec:solvers}.
The intention is to test both solvers of Chapter~\ref{sec:solvers} with both
discretizations of Chapter~\ref{sec:discretizations}.  However, the nonconforming
gradient-flow solver cannot be combined with the conforming discretization,
which cannot, without modifications, handle values not on the sphere.
We therefore end up with three combinations:
\begin{enumerate}
  \item The nonconforming gradient-flow solver of Chapter~\ref{sec:gradient-flow} for problems discretized with the nonconforming discretization
  (henceforth abbreviated as \textsf{Solve(GF)}/\textsf{Discr(N)}),
  \item the Riemannian trust-region solver of Chapter~\ref{sec:riemannianTR} with the nonconforming discretization (\textsf{Solve(TR)}/\textsf{Discr(N)}),
  \item the Riemannian trust-region solver with the conforming discretization
  (\textsf{\textsf{Solve(TR)}/Discr(C)}).
\end{enumerate}

We are interested in iteration numbers, wall times, and
how the nonconformity of the gradient-flow solver influences the final result.
In order to measure this nonconformity we use the quantity
\begin{equation*}
  \delta_1[\bv_h] \colonequals \int_{\Omega}  \cI_h^1 \big(\big| \abs{\bv_h}^2 - 1 \big|\big)\dv{x},
\end{equation*}
which is the quantity that appears in the theoretical constraint violation bound~\eqref{eq:constraint-estimate}.
As $\delta_1$ is zero for all iterates provided by the Riemannian trust-region method,
we show it only for the gradient-flow solver.

We again use the \CC/\textsc{Dune} implementation of the previous chapter.
For the wall-time measurements we used a standard laptop computer
with an AMD Ryzen~7 processor and 16\,GB of DDR4 RAM.
We set both solvers to iterate until the $H^1$-seminorm
of the correction drops below $\epsilon_\mathrm{stop} = 10^{-3}$.
Other than that, we used the same settings for the Riemannian trust-region solver
as before, i.e., initial trust-region radius $\Delta_0 = \frac{1}{2}$
and step acceptance parameters $\beta_1 = 0.9$ and $\beta_2 = 10^{-2}$.
For the gradient-flow solver we set the pseudo time step size $\tau$ to four times the maximum
element diameter $h$ of the grid.

\subsection{Inverse stereographic projection}

For the first test we again consider Problem~\ref{prob:inv_stereo},
where both domain $\Omega = (-\frac{1}{2},\frac{1}{2})^2$ and image space $S^2$ are two-dimensional.
This situation has a smooth solution in form of the inverse stereographic projection $\pi_\text{st}^{-1}$.

\subsubsection{Starting from close to the solution}

As a first test, we start the solvers from the nodal interpolation
of the function $\pi_\text{st}^{-1}$ itself.  This is the approach used for the discretization
error measurements of the previous chapter, only the termination criterion is now less strict.
Table~\ref{tab:inv_stereo-S} shows the iteration numbers.
Not surprisingly, convergence is very fast. Indeed, both solvers rarely need more than a single iteration
to reach a situation where the termination criterion holds.
Also, the constraint violation introduced by the gradient-flow solver remains very small.

As the iteration histories are so short we omit the wall-time measurements.

\begin{center}
\begin{table}\footnotesize
\addtolength{\tabcolsep}{-1.5pt}
\begin{tabular}{c r | c c  c c | c c | c  c }\footnotesize
  & &\multicolumn{4}{c|}{\textsf{Solve(GF)}/\textsf{Discr(N)}}    &\multicolumn{2}{c|}{\textsf{Solve(TR)}/\textsf{Discr(N)}} & \multicolumn{2}{c}{\textsf{Solve(TR)}/\textsf{Discr(C)}}\\
  & &\multicolumn{2}{c}{($p=1$)}&\multicolumn{2}{c|}{($p=2$)}    &($p = 1$)&($p=2$)    & ($p= 1$)&($p=2$)        \\
  $r$ &$\abs{\Th}$ &\#Iter &$\delta_\mathrm{1}[\bu_h]$ &\#Iter &$\delta_\mathrm{1}[\bu_h]$
  &\#Iter
  &\#Iter
  &\#Iter
  &\#Iter \\ \hline\hline
  1&  8         &1  & \num{2.776e-17} &3 & \num{2.776e-17}   & 1    &2  &1  &2    \\
  2&  32        &5  & \num{1.196e-06} &2 & \num{7.235e-09}   & 2    &2  &2  &2    \\
  3&  128       &4  & \num{4.370e-08} &1 & \num{1.551e-11}   & 2    &1  &2  &1    \\
  4&  512       &1  & \num{5.197e-10} &1 & \num{2.508e-14}   & 1    &1  &1  &1    \\
  5&  2\,048    &1  & \num{1.002e-11} &1 & \num{1.566e-16}   & 1    &1  &1  &1    \\
  6&  8\,192    &1  & \num{1.756e-13} &1 & \num{1.256e-16}   & 1    &1  &1  &1    \\
  7&  32\,768   &1  & \num{2.986e-15} &1 & \num{1.226e-16}   & 1    &1  &1  &1    \\
  8&  131\,072  &1  & \num{1.510e-16} &  &            & 1    &   &1  &   \\
  \hline
\end{tabular}%
\caption{Iteration numbers and unit-length constraint violation $\delta_1$
 for Problem~\ref{prob:inv_stereo}, starting close to the solution}
\label{tab:inv_stereo-S}
\end{table}
\end{center}

\subsubsection{Starting from further away}

To challenge the solvers a bit more, we now construct an initial iterate
that is further away from the solution.
For this, define the scalar perturbation function
\begin{equation}
  \widetilde{p}_n(x)
  \colonequals
  \cos(3\pi x_1) \cdot 4^{n} \prod\limits_{i=1}^{n} \Big( x_i^2 - \frac{1}{4} \Big),
  \label{eq:perturbationFunction}
\end{equation}
and use it to define the initial iterate
\begin{equation}
\label{eq:initial_iterate_stereographic}
 \bu^0_h
 \colonequals
 \mathcal{I}_h\Bigg(\frac{\Big(\pi_\text{st}^{-1} + \widetilde{p}_2 \cdot \begin{psmallmatrix} 1 \\ 0 \\ 0 \end{psmallmatrix}\Big)}
    {\abs[\Big]{\pi_\text{st}^{-1} + \widetilde{p}_2 \cdot \begin{psmallmatrix} 1 \\ 0 \\ 0 \end{psmallmatrix}}}\Bigg).
\end{equation}
See Figure~\ref{fig:radial-a-perturbed} for how this looks like.
The perturbation function~\eqref{eq:perturbationFunction} vanishes on the domain boundary $\partial\Omega$,
and therefore the initial iterate still satisfies the boundary condition
$\bu_h^0 = \pi_\text{st}^{-1}$ on $\partial \Omega$.
Also, the new initial iterate is still in the same homotopy class as $\pi_\text{st}^{-1}$.

    \begin{center}
      \begin{table}\footnotesize
      \addtolength{\tabcolsep}{-1.5pt}
      \begin{tabular}{c r | c c  c c | c c | c  c }\footnotesize
        & &\multicolumn{4}{c|}{\textsf{Solve(GF)}/\textsf{Discr(N)}}    &\multicolumn{2}{c|}{\textsf{Solve(TR)}/\textsf{Discr(N)}} & \multicolumn{2}{c}{\textsf{Solve(TR)}/\textsf{Discr(C)}}\\
        & &\multicolumn{2}{c}{($p=1$)}&\multicolumn{2}{c|}{($p=2$)}    &($p = 1$)&($p=2$)    & ($p= 1$)&($p=2$)        \\
        $r$ &$\abs{\Th}$ &\#Iter &$\delta_\mathrm{1}[\bu_h]$ &\#Iter &$\delta_\mathrm{1}[\bu_h]$
        &\#Iter
        &\#Iter
        &\#Iter
        &\#Iter \\ \hline\hline
        1&  8         & 14  & \num{3.630e-02} &  10   & \num{7.331e-02}      &  4  & 4  &  4    &  5  \\
        2&  32        & 15  & \num{4.485e-02} &  15   & \num{4.421e-02}      &  4  & 4  &  4    &  5  \\
        3&  128       & 24  & \num{1.837e-02} &  25   & \num{1.798e-02}      &  4  & 4  &  4    &  4  \\
        4&  512       & 43  & \num{1.032e-02} &  43   & \num{1.026e-02}      &  4  & 4  &  4    &  4  \\
        5&  2048      & 79  & \num{5.480e-03} &  79   & \num{5.472e-03}      &  4  & 4  &  4    &  4  \\
        6&  8192      &151  & \num{2.832e-03} & 151   & \num{2.831e-03}      &  4  & 4  &  4    &  4  \\
        7&  32\,768   &296  & \num{1.441e-03} & 296   & \num{1.440e-03}      &  4  & 4  &  4    &  4  \\
        8&  131\,072  &585  & \num{7.267e-04} &       &                &  4  &    &  4    &     \\
        \hline
      \end{tabular}%
      \caption{Iteration numbers and unit-length constraint violation $\delta_1$ for Problem~\ref{prob:inv_stereo}
      with initial iterate $\bu_0^h$ given by~\eqref{eq:initial_iterate_stereographic}}
      \label{tab:inv_stereo-S-perturbation}
      \end{table}
      \end{center}

    \begin{center} 
      \begin{table}\footnotesize
      \addtolength{\tabcolsep}{-1.5pt}
      \footnotesize
      \begin{tabular}{c r | c c | c c | c  c }\footnotesize
        & &\multicolumn{2}{c|}{\textsf{Solve(GF)}/\textsf{Discr(N)}}    &\multicolumn{2}{c|}{\textsf{Solve(TR)}/\textsf{Discr(N)}} & \multicolumn{2}{c}{\textsf{Solve(TR)}/\textsf{Discr(C)}}\\
        & &($p=1$)&($p=2$)    &($p = 1$)&($p=2$)    & ($p= 1$)&($p=2$)        \\
          $r$ &$\abs{\Th}$ & time [s] & time [s] & time [s] & time [s] & time [s] & time [s]\\ \hline\hline
          4&  512      &<1   &1     & \timenum{ 0.444366} & \timenum{ 0.565688} & \timenum{0.479742}   & \timenum{ 0.620856} \\
          5&  2\,048   &2    &8     & \timenum{ 0.68524 } & \timenum{ 1.21698 } & \timenum{0.819352}   & \timenum{ 1.38872 } \\
          6&  8\,192   &13   &99    & \timenum{ 1.64655 } & \timenum{ 3.40857 } & \timenum{1.99162 }   & \timenum{ 4.18479 } \\
          7&  32\,768  &152  &1269  & \timenum{ 5.38498 } & \timenum{12.4944  } & \timenum{6.53422 }   & \timenum{15.7255  } \\
          8&  131\,072 &1880 &      & \timenum{20.4656  } &                     & \timenum{25.159  }   &  \\
        \hline
      \end{tabular}%
      \caption{Wall-times for Problem~\ref{prob:inv_stereo} with perturbation,
      with initial iterate $\bu_0^h$ given by~\eqref{eq:initial_iterate_stereographic}}
       \label{tab:inv_stereo-S-times-perturbation}
      \end{table}
      \end{center}

Table~\ref{tab:inv_stereo-S-perturbation} shows the solver performance results. One can see that the trust-region method
still only needs a low number of iterations (4--5) to reach the required accuracy,
independent of the grid refinement and the approximation order.
The gradient-flow solver, on the other hand, needs much larger iteration numbers to reach the same accuracy.
This is caused by the step size choice $\tau = 4 h$, which leads to very small steps~$\tau$
for finer grids. Indeed, the iteration numbers seem to double from one grid to the next,
reflecting the linear dependence of $\tau$ on $h$.
The coupling of $\tau$ and $h$, however, is necessary in order to keep a bound on
the constraint violation $\delta_1$ according to the estimate~\eqref{eq:constraint-estimate}.

Table~\ref{tab:inv_stereo-S-perturbation} shows that bounding the constraint violation
in this way does work in practice.
The constraint violations of the minimizers computed by the gradient-flow method
remain in the range of $10^{-4}$ to $10^{-3}$, which will be negligible for many practical purposes.
In fact, as predicted by~\eqref{eq:constraint-estimate}, the constraint violation is proportional to the time-step size:
As we have coupled $\tau$ to be proportional to $h$, the violation is roughly reduced by a factor of $2$ for each grid refinement.

High iteration numbers would be unproblematic if the gradient-flow iterations were cheap.
However, both solvers have roughly the same cost per iteration, dominated by
having to solve a linear system of equations for each iteration.
The matrix of the trust-region method is symmetric and positive definite,
but it is iteration-dependent, and therefore
has to be reassembled at each iteration. The linear system of the gradient-flow method
is a saddle-point problem, on the other hand, because (in our implementation)
the tangentiality of the correction is enforced via Lagrange multipliers.
The benefit of this is that the two diagonal blocks of the matrix are independent
of the iteration, and have to be assembled only once. On the downside, the problem
has more unknowns than the formulation in local coordinates of the tangent space.
An alternative implementation could formulate the tangent problem of the
gradient-flow method in local coordinates of the tangent space, which would
get rid of the Lagrange multipliers. Then, however, the entire matrix
would have to be reassembled at each step as well.

The net effect of iteration numbers and time-per iteration can be seen in Table~\ref{tab:inv_stereo-S-times-perturbation}.
The wall-time of the trust-region solver scales roughly linearly with the number of
degrees of freedom. This is the combination of the fact that the outer trust-region solver
needs a resolution-independent number of iterations, and that the inner solver is a multigrid
solver with optimal complexity. The gradient-flow solver does not have these features,
and therefore the wall-time it requires increases much faster as the grid gets finer.
The difference between the two is in the range of two orders of magnitude on the finer grids.

\subsection{Radial projection}

In the next sequence of tests we consider the benchmark problems
\ref{prob:singular}\ref{subprob:radial33} and \ref{prob:singular}\ref{subprob:radial22},
which involve harmonic maps with one singularity.
We will see that this singularity influences the solver behavior.

\subsubsection{The case $n=m=3$}
\label{sec:solver_radial_3}
We start with Problem~\ref{prob:singular}\ref{subprob:radial33}, which asks for a harmonic map on $\Omega = (-\frac{1}{2}, \frac{1}{2})^3$
with image in $S^2$, with boundary data given by the radial projection $\bu_{\odot}: x\mapsto \frac{x}{\abs{x}}$.
This projection is also a solution of the problem. It is singular at the origin,
but nevertheless an element of $H^1(\Omega;S^2)$.

As the first test we start the solvers directly from~\eqref{eq:radial_projection_with_origin_tweak_3d},
which is essentially the radial projection itself.
Table~\ref{tab:radial33-S} shows the number of iterations.
One can see that this problem is more difficult than the previous one:
The number of required trust-region iterations is still bounded independent
from the grid resolution, but $5$ to $8$ iterations are now needed to reach
the termination criterion even for this good initial iterate.
The gradient-flow solver, on the other hand, quickly requires three-digit
iteration numbers, and the inverse proportional dependence on the
grid element diameter~$h$ can be observed again.

Table~\ref{tab:radial33-S} also shows the constraint violation $\delta_1$ of the gradient-flow solver. As in the previous example it decreases with each refinement step. Unlike previously,
the reduction is even better than what is expected from the bound~\eqref{eq:constraint-estimate}.

The wall-time measurements in Table~\ref{tab:radial33-Times} reflect
the different iteration numbers. Even considering that each refinement now
multiplies the number of degrees of freedom by~$8$,
the time complexity of the trust-region is not quite optimal anymore.
Still, it is between one and two orders of magnitude faster than the gradient-flow solver.

\begin{center}
  \begin{table}\footnotesize
  \addtolength{\tabcolsep}{-1.5pt}
  \begin{tabular}{c r | c c  c c | c c | c  c }\footnotesize
    & &\multicolumn{4}{c|}{\textsf{Solve(GF)}/\textsf{Discr(N)}}    &\multicolumn{2}{c|}{\textsf{Solve(TR)}/\textsf{Discr(N)}} & \multicolumn{2}{c}{\textsf{Solve(TR)}/\textsf{Discr(C)}}\\
    & &\multicolumn{2}{c}{($p=1$)}&\multicolumn{2}{c|}{($p=2$)}    &($p = 1$)&($p=2$)    & ($p= 1$)&($p=2$)        \\
    $r$ &$\abs{\Th}$ &\#Iter &$\delta_\mathrm{1}[\bu_h]$ &\#Iter &$\delta_\mathrm{1}[\bu_h]$
    &\#Iter
    &\#Iter
    &\#Iter
    &\#Iter \\ \hline\hline
    1&  48            &  1  & \num{4.163e-17} & 58  & \num{2.015e-03}     &1   & 4  & 6  &5    \\
    2&  384           &353  & \num{1.171e-03} & 90  & \num{4.888e-04}     &5   & 8  & 6  &5    \\
    3&  3072          &244  & \num{2.832e-04} &159  & \num{6.022e-05}     &5   & 8  & 7  &5    \\
    4&  24\,576       &387  & \num{4.222e-05} &299  & \num{6.850e-06}     &6   & 7  & 7  &5    \\
    5&  196\,608      &750  & \num{5.323e-06} &     &                     &6   & 7  & 7  &5    \\
    \hline
  \end{tabular}%
  \caption{Iteration numbers and unit-length constraint violation for Problem~\ref{prob:singular}\ref{subprob:radial33},
  starting from close to the solution}
  \label{tab:radial33-S}
  \end{table}
  \end{center}

    \begin{center}
      \begin{table}\footnotesize
      \addtolength{\tabcolsep}{-1.5pt}
      \footnotesize
      \begin{tabular}{c r | c c | c c | c  c }\footnotesize
        & &\multicolumn{2}{c|}{\textsf{Solve(GF)}/\textsf{Discr(N)}}    &\multicolumn{2}{c|}{\textsf{Solve(TR)}/\textsf{Discr(N)}} & \multicolumn{2}{c}{/\textsf{Solve(TR)}/\textsf{Discr(C)}}\\
        & &($p=1$)&($p=2$)    &($p = 1$)&($p=2$)    & ($p= 1$)&($p=2$)        \\
          $r$ &$\abs{\Th}$ & time [s] & time [s] & time [s] & time [s] & time [s] & time [s]\\ \hline\hline
          3&  3072     & \timenum{1.000e+01} & \timenum{1.660e+02}  & \timenum{ 1.11307} & \timenum{  7.29428} & \timenum{ 1.40861} & \timenum{  5.48709}   \\
          4&  24\,576  & \timenum{1.770e+02} & \timenum{3.358e+03}  & \timenum{ 6.65815} & \timenum{ 69.9041 } & \timenum{ 9.48989} & \timenum{ 53.4734 }   \\
          5&  196\,608 & \timenum{3.509e+03} &                      & \timenum{60.6154 } & \timenum{776.859  } & \timenum{83.5153 } & \timenum{649.919  }   \\
        \hline
      \end{tabular}%
      \caption{Wall-times for Problem~\ref{prob:singular}\ref{subprob:radial33},
      starting from close to the solution}
      \label{tab:radial33-Times}
      \end{table}
      \end{center}

\subsubsection{Starting farther away from the solution} \label{ex:random}

As for the previous problem we now start the solver at
initial data that is further away from the discrete solution.
To construct an initial iterate, we reuse the scalar perturbation function
$\widetilde{p}_n$ defined in~\eqref{eq:perturbationFunction},
and define
\begin{equation}
 \bu^0_h(x)
 \colonequals
 \mathcal{I}_h
 \Bigg(
 \frac{x + \widetilde{p_3}(x)\begin{psmallmatrix}1 \\1\\0\end{psmallmatrix}}
 {\abs[\Big]{x + \widetilde{p_3}(x)\begin{psmallmatrix}1 \\1\\0\end{psmallmatrix}}}
 \Bigg).
  \label{eq:radial33-perturbedInitialIterate}
\end{equation}
This new initial iterate is visualized in Figure~\ref{fig:radial-b-perturbed}.
The singularity has moved to about $(0.148, 0.148, 0)^T$,
but by the construction of $\widetilde{p}_n$ in~\eqref{eq:perturbationFunction}, the map $\bu^0_h$ still satisfies the
boundary conditions $\bu^0_h = \bu_\odot$ on $\partial \Omega$.

\begin{figure}
  \centering
  \begin{subfigure}[b]{0.32\textwidth}
    \centering
    \includegraphics[width=0.75\textwidth]{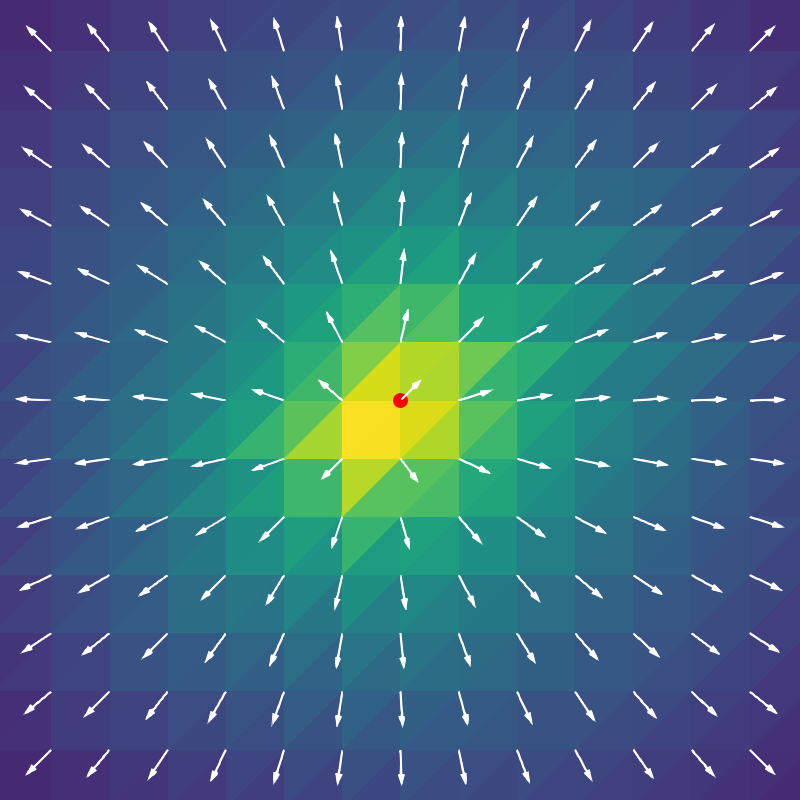}
    \caption{no perturbation, $r=4$}
  \end{subfigure}
  \begin{subfigure}[b]{0.32\textwidth}
    \centering
    \includegraphics[width=0.75\textwidth]{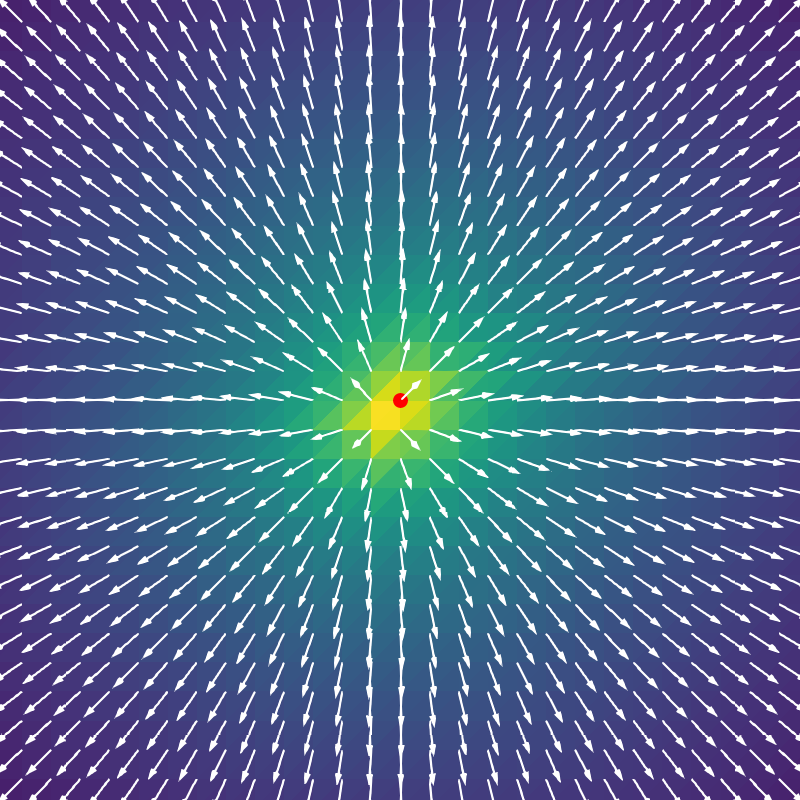}
    \caption{no perturbation, $r=5$}
  \end{subfigure}
  \begin{subfigure}[b]{0.32\textwidth}
    \centering
    \includegraphics[width=0.75\textwidth]{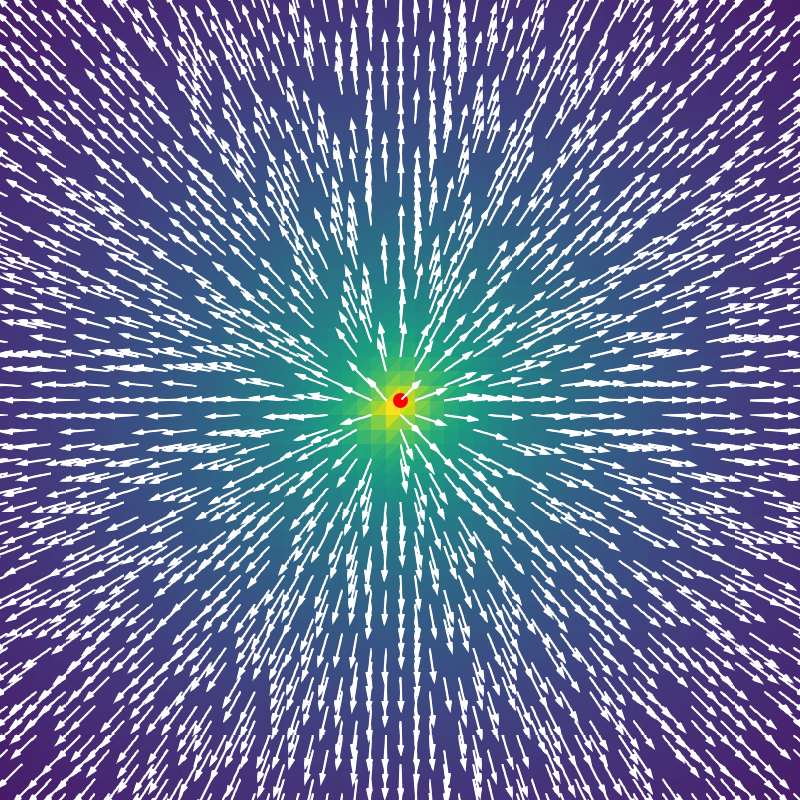}
    \caption{no perturbation, $r=6$}
  \end{subfigure}
  \begin{subfigure}[b]{0.32\textwidth}
    \centering
    \includegraphics[width=0.75\textwidth]{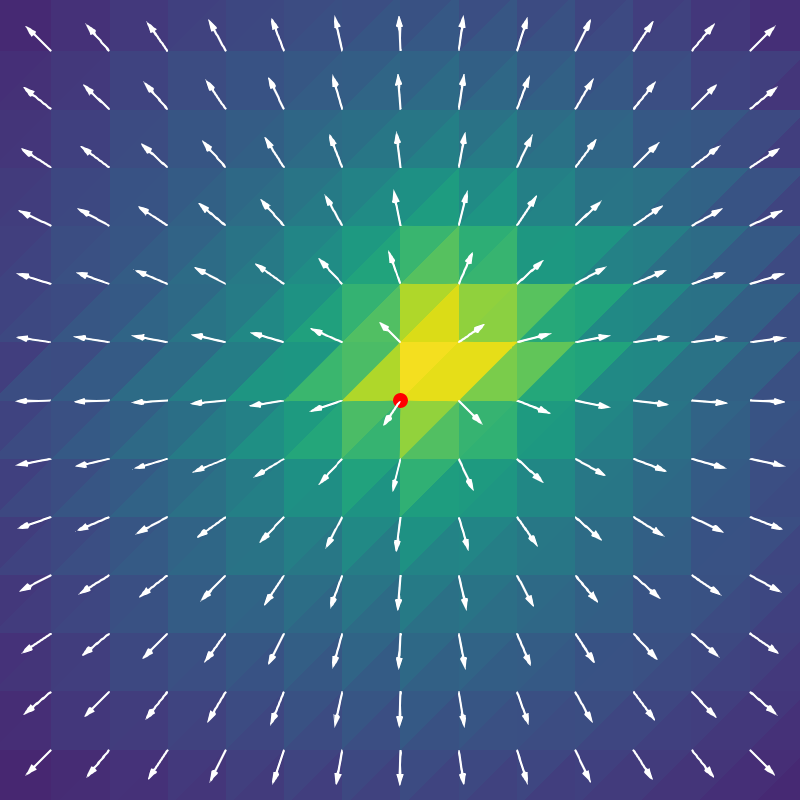}
    \caption{with perturbation, $r=4$}
  \end{subfigure}
  \begin{subfigure}[b]{0.32\textwidth}
    \centering
    \includegraphics[width=0.75\textwidth]{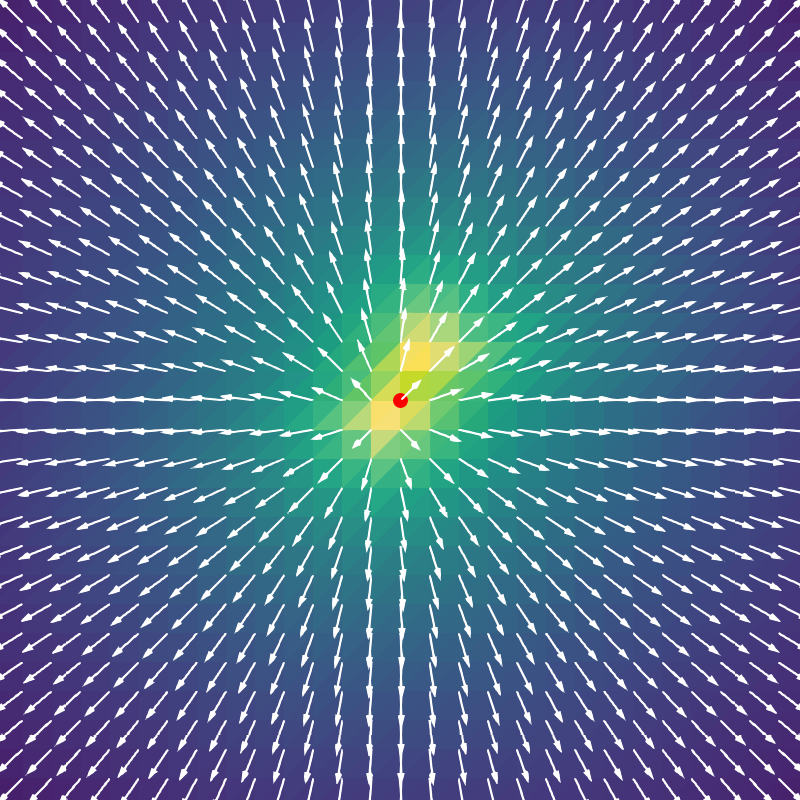}
    \caption{with perturbation, $r=5$}
  \end{subfigure}
  \begin{subfigure}[b]{0.32\textwidth}
    \centering
    \includegraphics[width=0.75\textwidth]{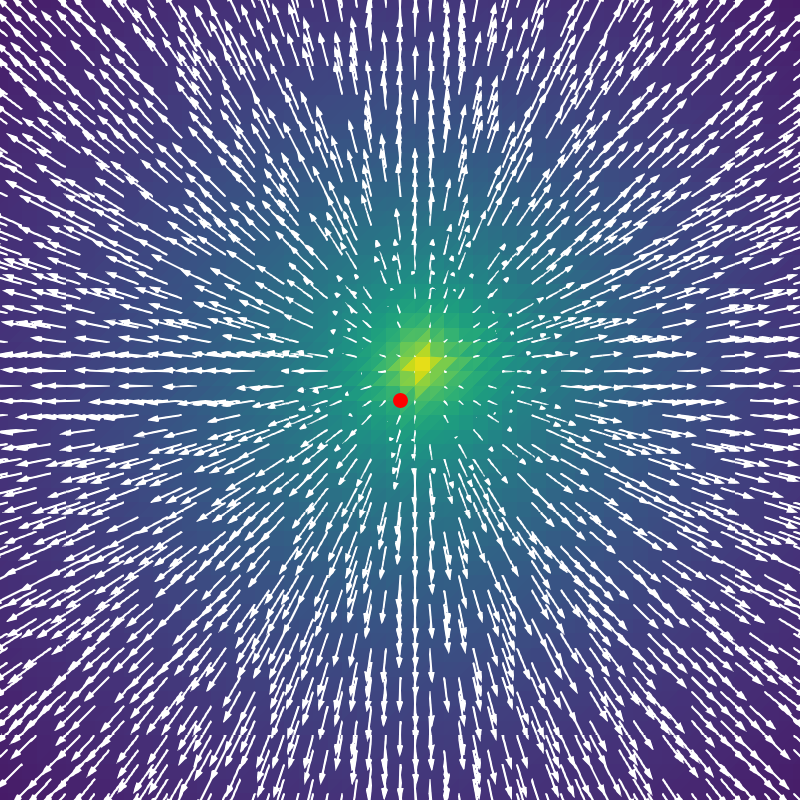}
    \caption{with perturbation, $r=6$}
  \end{subfigure}
  \caption{Nonconforming discretization: Limit configurations $\bu_h$
  for Problem~\ref{prob:singular}\ref{subprob:radial33} for initial
  iterate~\eqref{eq:radial_projection_with_origin_tweak_3d} close to the solution (top row),
  and with initial iterate \eqref{eq:radial33-perturbedInitialIterate} (bottom row)}
  \label{fig:radial33-nonconforming-initialIterateDependence}
\end{figure}

\begin{figure}
  \centering
  \begin{subfigure}[b]{0.32\textwidth}
    \centering
    \includegraphics[width=0.75\textwidth]{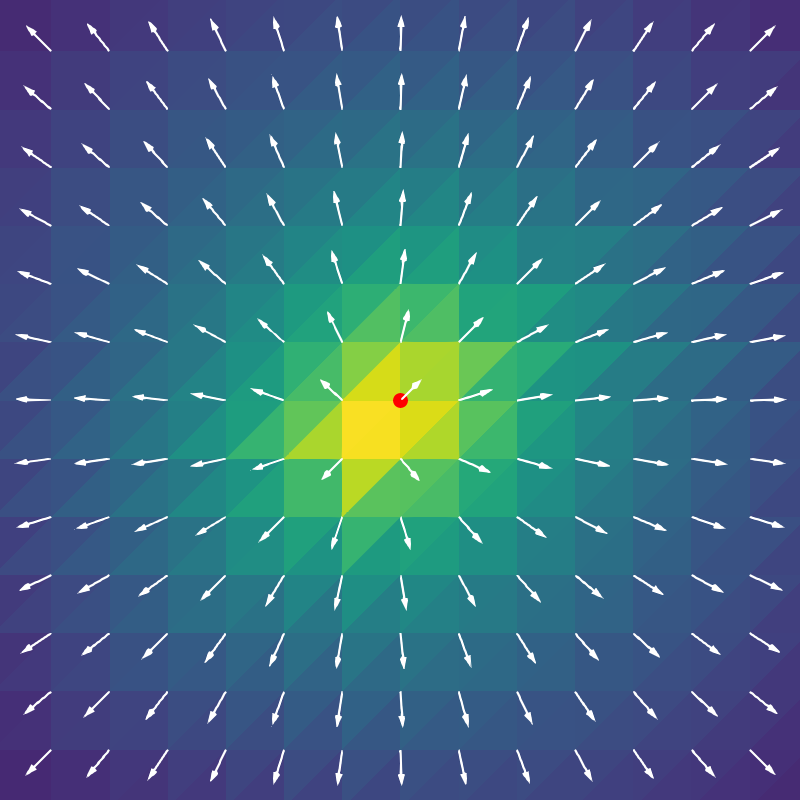}
    \caption{no perturbation, $r=4$}
  \end{subfigure}
  \begin{subfigure}[b]{0.32\textwidth}
    \centering
    \includegraphics[width=0.75\textwidth]{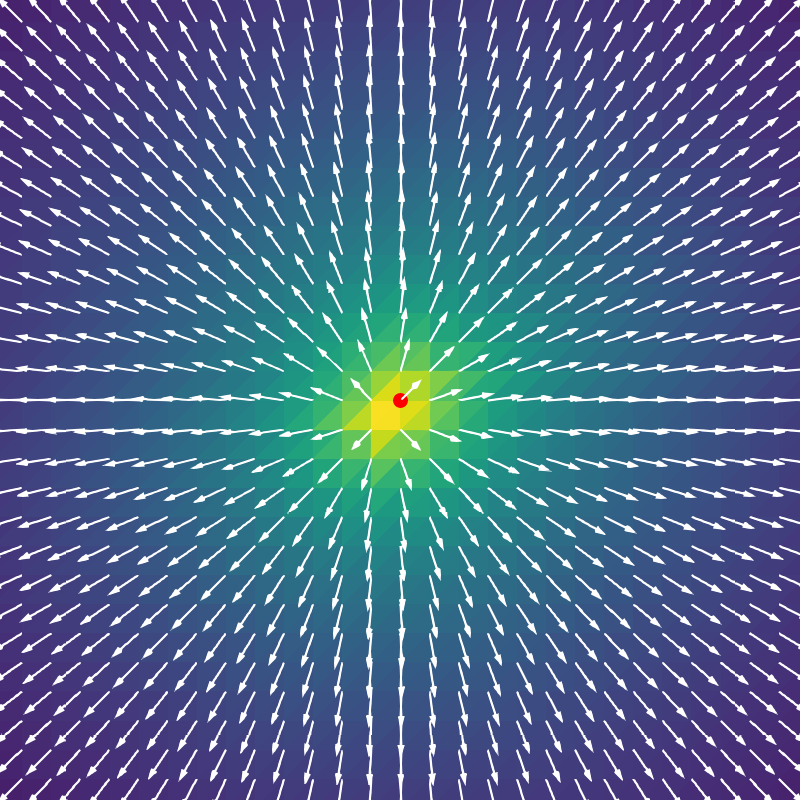}
    \caption{no perturbation, $r=5$}
  \end{subfigure}
  \begin{subfigure}[b]{0.32\textwidth}
    \centering
    \includegraphics[width=0.75\textwidth]{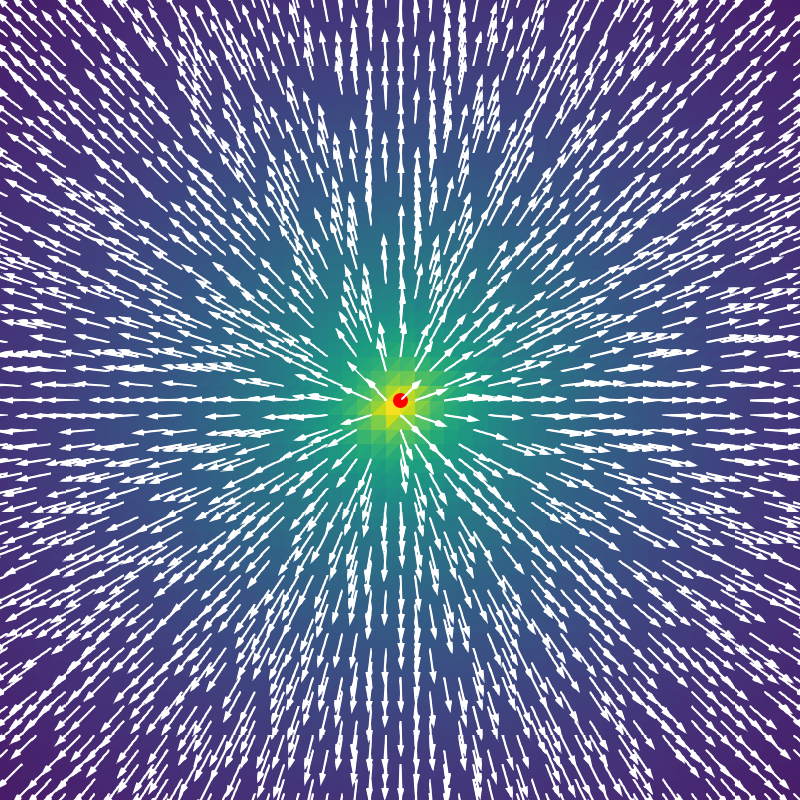}
    \caption{no perturbation, $r=6$}
  \end{subfigure}
  \begin{subfigure}[b]{0.32\textwidth}
    \centering
    \includegraphics[width=0.75\textwidth]{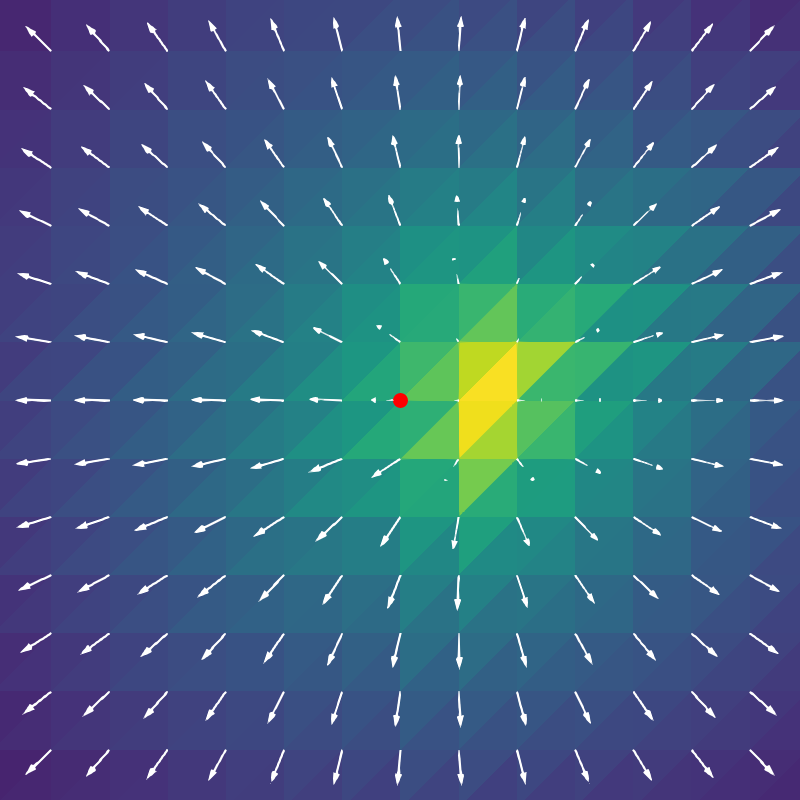}
    \caption{with perturbation, $r=4$}
  \end{subfigure}
  \begin{subfigure}[b]{0.32\textwidth}
    \centering
    \includegraphics[width=0.75\textwidth]{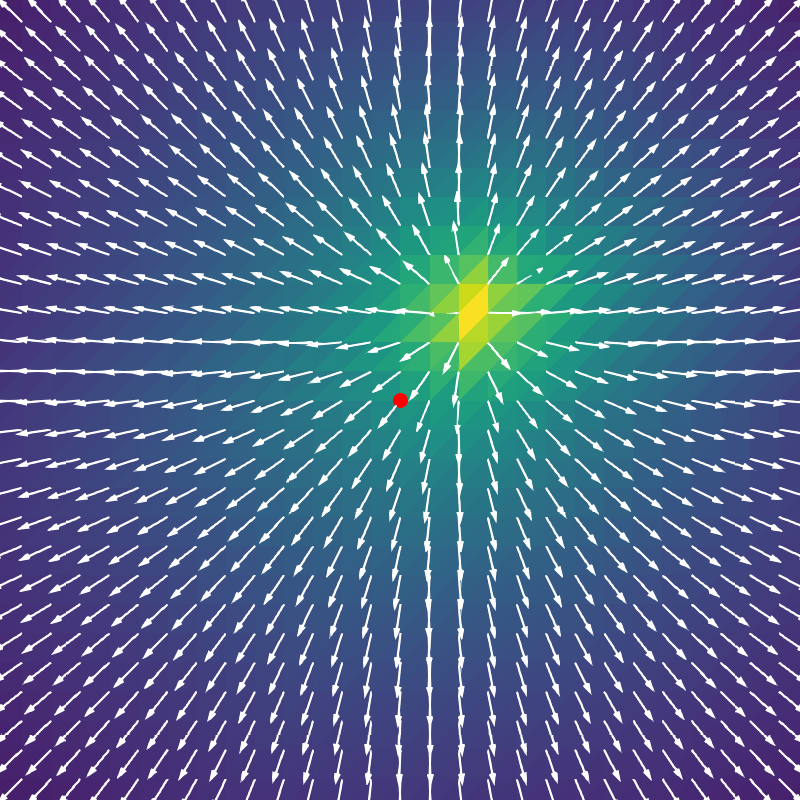}
    \caption{with perturbation, $r=5$}
  \end{subfigure}
  \begin{subfigure}[b]{0.32\textwidth}
    \centering
    \includegraphics[width=0.75\textwidth]{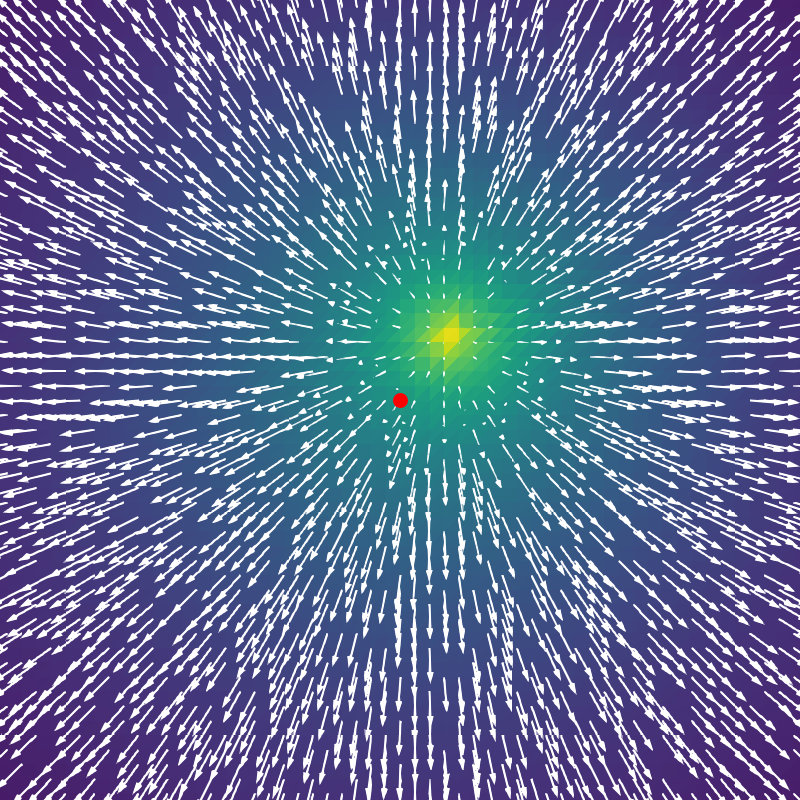}
    \caption{with perturbation, $r=6$}
  \end{subfigure}
  \caption{Conforming discretization: Limit configurations $\bu_h$
  for Problem~\ref{prob:singular}\ref{subprob:radial33} for initial iterate~\eqref{eq:radial_projection_with_origin_tweak_3d}
  close to the solution (top row), and with initial iterate
  \eqref{eq:radial33-perturbedInitialIterate} (bottom row)}
  \label{fig:radial33-conforming-initialIterateDependence}
\end{figure}

Table~\ref{tab:radial33-S-perturbation} shows the iteration numbers and nonconformity,
and Table~\ref{tab:radial33-Times-perturbation} shows the wall times.
The iteration numbers and wall times have increased further, but qualitatively
the behavior is still the same as for the previous examples.
However, looking at the computed configurations reveals a problem:
Figures~\ref{fig:radial33-nonconforming-initialIterateDependence}
and~\ref{fig:radial33-conforming-initialIterateDependence} show the limit configurations
on three different grids for the nonconforming and the conforming discretization,
respectively. In all pictures, a red dot marks the origin, i.e., the place where
the singularity should be. One can see that its final computed position is not
at the origin when starting from~\eqref{eq:radial33-perturbedInitialIterate}.
Rather, it seems to be stuck at where it was in $\bu_h^0$. This is not caused by
the relatively weak termination criterion---the results
do not change when setting $\epsilon_\text{stop}$ to $10^{-15}$. The effect is markedly
stronger for the conforming discretization, but it is also visible for the
nonconforming one. There seems to be an effect that obstructs the movement
of singularities. One conjecture is that the fact that singularities are restricted
to certain parts of an element in a projection-based finite element computation
introduces these obstructions, but we have no direct justification for this yet.

Our results are in line with~\cite{huWinther2009saddle}, where,
in a two-dimensional example, the authors
also observed different limit configurations depending on the initial iterate.

    \begin{table}\footnotesize
    \addtolength{\tabcolsep}{-1.5pt}
    \begin{center}
    \begin{tabular}{c r | c c  c c | c c | c  c }\footnotesize
      & &\multicolumn{4}{c|}{\textsf{Solve(GF)}/\textsf{Discr(N)}}    &\multicolumn{2}{c|}{\textsf{Solve(TR)}/\textsf{Discr(N)}} & \multicolumn{2}{c}{\textsf{Solve(TR)}/\textsf{Discr(C)}}\\
      & &\multicolumn{2}{c}{($p=1$)}&\multicolumn{2}{c|}{($p=2$)}    &($p = 1$)&($p=2$)    & ($p= 1$)&($p=2$)        \\
      $r$ &$\abs{\Th}$ &\#Iter &$\delta_\mathrm{1}[\bu_h]$ &\#Iter &$\delta_\mathrm{1}[\bu_h]$
      &\#Iter
      &\#Iter
      &\#Iter
      &\#Iter \\ \hline\hline
      1&  48            &   1  & \num{4.626e-17} & 115   & \num{2.025e-02}  &  1 &  9  &  6  &  6   \\
      2&  384           & 237  & \num{1.625e-02} &  79   & \num{2.287e-02}  & 10 &  8  & 10  & 14   \\
      3&  3072          & 125  & \num{7.342e-03} & 361   & \num{7.507e-03}  &  8 & 10  & 10  & 13   \\
      4&  24\,576       & 920  & \num{3.894e-03} &1345   & \num{3.897e-03}  & 11 & 13  &  9  &  8   \\
      5&  196\,608      &5077  & \num{2.009e-03} &       &                  & 10 & 20  & 11  & 10   \\
      \hline
    \end{tabular}%
    \end{center}
    \caption{Iteration numbers and unit-length violation for Problem~\ref{prob:singular}\ref{subprob:radial33} with initial iterate~\eqref{eq:radial33-perturbedInitialIterate}}
    \label{tab:radial33-S-perturbation}
    \end{table}

    \begin{table}\footnotesize
    \addtolength{\tabcolsep}{-1.5pt}
    \footnotesize
    \begin{center}
    \begin{tabular}{c r | c c | c c | c  c }\footnotesize
      & &\multicolumn{2}{c|}{\textsf{Solve(GF)}/\textsf{Discr(N)}}    &\multicolumn{2}{c|}{\textsf{Solve(TR)}/\textsf{Discr(N)}} & \multicolumn{2}{c}{\textsf{Solve(TR)}/\textsf{Discr(C)}}\\
      & &($p=1$)&($p=2$)    &($p = 1$)&($p=2$)    & ($p= 1$)&($p=2$)        \\
        $r$ &$\abs{\Th}$ & time [s] & time [s] & time [s] & time [s] & time [s] & time [s]\\ \hline\hline
        1&  48          & <1      &3     & \timenum{  0.780034} & \timenum{   1.45382} & \timenum{  0.212625} & \timenum{   0.454979} \\
        2&  384         & 6       &7     & \timenum{  1.065   } & \timenum{   1.9858 } & \timenum{  0.587439} & \timenum{   3.27731 } \\
        3&  3072        & 6       &364   & \timenum{  2.55349 } & \timenum{  11.7205 } & \timenum{  2.52222 } & \timenum{  14.5316  } \\
        4&  24\,576     & 413     &15040 & \timenum{ 14.6909  } & \timenum{ 144.848  } & \timenum{ 17.0718  } & \timenum{ 101.011   } \\
        5&  196\,608    & 23\,700 &      & \timenum{105.062   } & \timenum{2201.44   } & \timenum{141.482   } & \timenum{1326.09    } \\
      \hline
    \end{tabular}%
    \caption{Wall-times for Problem~\ref{prob:singular}\ref{subprob:radial33}
             with initial iterate~\eqref{eq:radial33-perturbedInitialIterate}}
    \label{tab:radial33-Times-perturbation}
    \end{center}
    \end{table}

\subsubsection{The case $n=m=2$}

Next, we investigate Problem~\ref{prob:singular}\ref{subprob:radial22}, the approximation of a harmonic map
from $(-\frac{1}{2}, \frac{1}{2})^2$ to $S^1$ with boundary data given by the radial projection
$\bu_\odot$ now on this two-dimensional domain.  Since the solution of this,
i.e., $\bu_\odot$ itself, is now even less regular than in the three-dimensional scenario, we expect
additional numerical difficulties.  However, the numerical results suggest that these fears
are unfounded.

We again start from the nodal
interpolation of the solution $x \mapsto \frac{x}{\abs{x}}$, with the modification
described in Chapter~\ref{ex:radial22}.
Table~\ref{tab:radial22-S} shows the total number of iterations as well as the unit-length constraint violation $\delta_1$.
The number of iterations of both methods is comparable to the three-dimensional situation
of Chapter~\ref{sec:solver_radial_3}, which had an $H^1$-solution (Table~\ref{tab:radial33-S}).
However, the nonconforming discretization seems to slightly reduce the number of total iteration steps of the Riemannian trust-region method.

  \begin{table}\footnotesize
  \addtolength{\tabcolsep}{-1.5pt}
  \begin{center}
  \begin{tabular}{c r | c c  c c | c c | c  c }\footnotesize
    & &\multicolumn{4}{c|}{\textsf{Solve(GF)}/\textsf{Discr(N)}}    &\multicolumn{2}{c|}{\textsf{Solve(TR)}/\textsf{Discr(N)}} & \multicolumn{2}{c}{\textsf{Solve(TR)}/\textsf{Discr(C)}}\\
    & &\multicolumn{2}{c}{($p=1$)}&\multicolumn{2}{c|}{($p=2$)}    &($p = 1$)&($p=2$)    & ($p= 1$)&($p=2$)        \\
    $r$ &$\abs{\Th}$ &\#Iter &$\delta_\mathrm{1}[\bu_h]$ &\#Iter &$\delta_\mathrm{1}[\bu_h]$
    &\#Iter
    &\#Iter
    &\#Iter
    &\#Iter \\ \hline\hline
    1&  8         &    1 & \num{5.551e-17} & 12  & \num{5.551e-17}   &1     &5  &4  &5    \\
    2&  32        &   18 & \num{4.752e-03} & 21  & \num{4.611e-04}   &5     &6  &7  &7    \\
    3&  128       &   33 & \num{1.507e-03} & 36  & \num{1.429e-04}   &5     &5  &9  &7    \\
    4&  512       &   61 & \num{3.042e-04} & 64  & \num{3.003e-05}   &6     &5  &9  &7    \\
    5&  2048      &  117 & \num{5.116e-05} &118  & \num{5.238e-06}   &6     &5  &9  &7    \\
    6&  8192      &  226 & \num{7.817e-06} &228  & \num{8.241e-07}   &6     &5  &9  &7    \\
    7&  32\,768   &  445 & \num{1.130e-06} &446  & \num{1.220e-07}   &6     &6  &9  &7    \\
    8&  131\,072  &  883 & \num{1.581e-07} &     &                   &6     &   &9  &   \\
    \hline
  \end{tabular}%
  \end{center}
  \caption{Iteration numbers and unit-length violation for Problem~\ref{prob:singular}\ref{subprob:radial22},
  starting close to the solution}
  \label{tab:radial22-S}
  \end{table}

    \begin{table}\footnotesize
    \addtolength{\tabcolsep}{-1.5pt}
    \footnotesize
    \begin{center}
    \begin{tabular}{c r | c c | c c | c  c }\footnotesize
      & &\multicolumn{2}{c|}{\textsf{Solve(GF)}/\textsf{Discr(N)}}    &\multicolumn{2}{c|}{\textsf{Solve(TR)}/\textsf{Discr(N)}} & \multicolumn{2}{c}{\textsf{Solve(TR)}/\textsf{Discr(C)}}\\
      & &($p=1$)&($p=2$)    &($p = 1$)&($p=2$)    & ($p= 1$)&($p=2$)        \\
        $r$ &$\abs{\Th}$ & time [s] & time [s] & time [s] & time [s] & time [s] & time [s]\\ \hline\hline
        4&  512      & <1                  & \timenum{1}        & \timenum{ 0.724762} & \timenum{ 1.40011} & \timenum{ 1.37069} & \timenum{ 2.48307} \\
        5&  2048     & \timenum{1}         & \timenum{6}        & \timenum{ 1.50181 } & \timenum{ 4.11088} & \timenum{ 2.70327} & \timenum{ 6.49168} \\
        6&  8192     & \timenum{9}         & \timenum{7.6e+01}  & \timenum{ 4.47555 } & \timenum{13.2629 } & \timenum{ 7.36816} & \timenum{18.8744 } \\
        7&  32\,768  & \timenum{9.5e+01}   & \timenum{8.58e+02} & \timenum{15.06    } & \timenum{51.0858 } & \timenum{22.4757 } & \timenum{60.8347 } \\
        8&  131\,072 & \timenum{1.244e+03} &                    & \timenum{48.2341  } &                    & \timenum{71.8069 } &  \\
      \hline
    \end{tabular}%
    \end{center}
    \caption{Wall-times for Problem~\ref{prob:singular}\ref{subprob:radial22},
    starting close to the solution}
    \label{tab:radial22-times}
    \end{table}

The wall times, shown in Table~\ref{tab:radial22-times},
match the results of the earlier test problems: While the discrete gradient flow
solver is faster for smaller problems, the trust-region solver outperforms it
by orders of magnitude once the grid reaches a certain size.

\subsubsection{$n=m=2$: Starting farther away from the solution} \label{ex:2b-pertubed}

Similar to the previous examples we use a perturbation of $\bu_\odot$ as a starting configuration (see Figure~\ref{fig:radial-c})
\begin{equation*}
 \bu^0_h(x)
 \colonequals
 \mathcal{I}_h \Bigg(
  \frac{x + \widetilde{p_2}(x)\begin{psmallmatrix} 1 \\ 0\end{psmallmatrix}}
  {\abs[\big]{x + \widetilde{p_2}(x)\begin{psmallmatrix} 1 \\ 0\end{psmallmatrix}}} \Bigg).
\end{equation*}

        \begin{table}\footnotesize
        \addtolength{\tabcolsep}{-1.5pt}
        \begin{center}
        \begin{tabular}{c r | c c  c c | c c | c  c }\footnotesize
          & &\multicolumn{4}{c|}{\textsf{Solve(GF)}/\textsf{Discr(N)}}  &\multicolumn{2}{c|}{\textsf{Solve(TR)}/\textsf{Discr(N)}} & \multicolumn{2}{c}{\textsf{Solve(TR)}/\textsf{Discr(C)}}\\
          & &\multicolumn{2}{c}{($p=1$)}&\multicolumn{2}{c|}{($p=2$)}    &($p = 1$)&($p=2$)    & ($p= 1$)&($p=2$)        \\
          $r$ &$\abs{\Th}$ &\#Iter &$\delta_\mathrm{1}[\bu_h]$ &\#Iter &$\delta_\mathrm{1}[\bu_h]$
          &\#Iter
          &\#Iter
          &\#Iter
          &\#Iter \\ \hline\hline
          1&  8         &   1 & \num{5.551e-17} &  32  & \num{1.284e-02}   & 1  & 5  & 4   & 7\\
          2&  32        &  66 & \num{2.204e-02} &  42  & \num{2.346e-02}   & 5  & 6  & 6   &10\\
          3&  128       &  68 & \num{9.836e-03} &  95  & \num{9.501e-03}   & 7  & 7  & 5   & 6\\
          4&  512       & 157 & \num{5.705e-03} & 215  & \num{5.763e-03}   & 7  & 7  & 9   &11\\
          5&  2048      & 378 & \num{3.105e-03} & 492  & \num{3.199e-03}   & 7  & 9  & 8   & 9\\
          6&  8192      & 903 & \num{1.661e-03} &1130  & \num{1.704e-03}   &10  & 9  & 7   & 9\\
          7&  32\,768   &2132 & \num{8.669e-04} &2591  & \num{8.817e-04}   & 9  & 9  & 1   & 9\\
          8&  131\,072  &4965 & \num{4.443e-04} &      &                   & 8  &    & 9   &  \\
          \hline
        \end{tabular}%
        \end{center}
        \caption{Iteration numbers and unit-length violation for Problem~\ref{prob:singular}\ref{subprob:radial22} starting away from the solution}
        \label{tab:radial22-S-perturbed}
        \end{table}

        \begin{table}\footnotesize
        \addtolength{\tabcolsep}{-1.5pt}
        \footnotesize
        \begin{center}
        \begin{tabular}{c r | c c | c c | c  c }\footnotesize
          & &\multicolumn{2}{c|}{\textsf{Solve(GF)}/\textsf{Discr(N)}}    &\multicolumn{2}{c|}{\textsf{Solve(TR)}/\textsf{Discr(N)}} & \multicolumn{2}{c}{\textsf{Solve(TR)}/\textsf{Discr(C)}}\\
          & &($p=1$)&($p=2$)    &($p = 1$)&($p=2$)    & ($p= 1$)&($p=2$)        \\
            $r$ &$\abs{\Th}$ & time [s] & time [s] & time [s] & time [s] & time [s] & time [s]\\ \hline\hline
            4&  512      &<1    &2     & \timenum{ 0.974793} & \timenum{ 2.16622 } & \timenum{ 1.33396 } & \timenum{ 3.2145  } \\
            5&  2048     &3     &25    & \timenum{ 2.14577 } & \timenum{ 7.62299 } & \timenum{ 2.40497 } & \timenum{ 7.71195 } \\
            6&  8192     &38    &368   & \timenum{ 7.2623  } & \timenum{23.6048  } & \timenum{ 6.06336 } & \timenum{24.1517  } \\
            7&  32\,768  &450   &4816  & \timenum{22.6324  } & \timenum{79.2613  } & \timenum{27.0292  } & \timenum{81.2429  } \\
            8&  131\,072 &6941  &      & \timenum{67.5566  } &           &77.2647    &           \\
          \hline
        \end{tabular}%
        \end{center}
        \caption{Wall-times for Problem~\ref{prob:singular}\ref{subprob:radial22},
        starting away from the solution}
        \label{tab:radial22-times-perturbed}
        \end{table}

Very little changes for the trust-region solver: The iteration numbers increase only a bit,
and the wall time is increased by about 30\,\%. For the gradient-flow solver, the increase
of iteration numbers is above 500\,
with a corresponding increase of wall time. Finally, one can see from Table~\ref{tab:radial22-S-perturbed}
that the constraint violation accumulated by the gradient-flow solver is now
a factor 1000 larger than when starting from the radial projection. In absolute numbers
it is still small, though.

Curiously, in this two-dimensional setting we could not reproduce the effect shown
in Figure~\ref{fig:radial33-conforming-initialIterateDependence}, where the singularity would get stuck
away from its optimal position, regardless of the grid resolution and solver precision.
In the two-dimensional experiments, the limit position of the singularity produced
by both solvers and discretizations would always be right at the origin,
where it should be.

\subsection{Harmonic maps with multiple singularities}
\label{sec:solver_multiple_singularities}

\begin{figure}
  \centering
  \begin{subfigure}[b]{\textwidth}
    \includegraphics[width=0.25\textwidth]{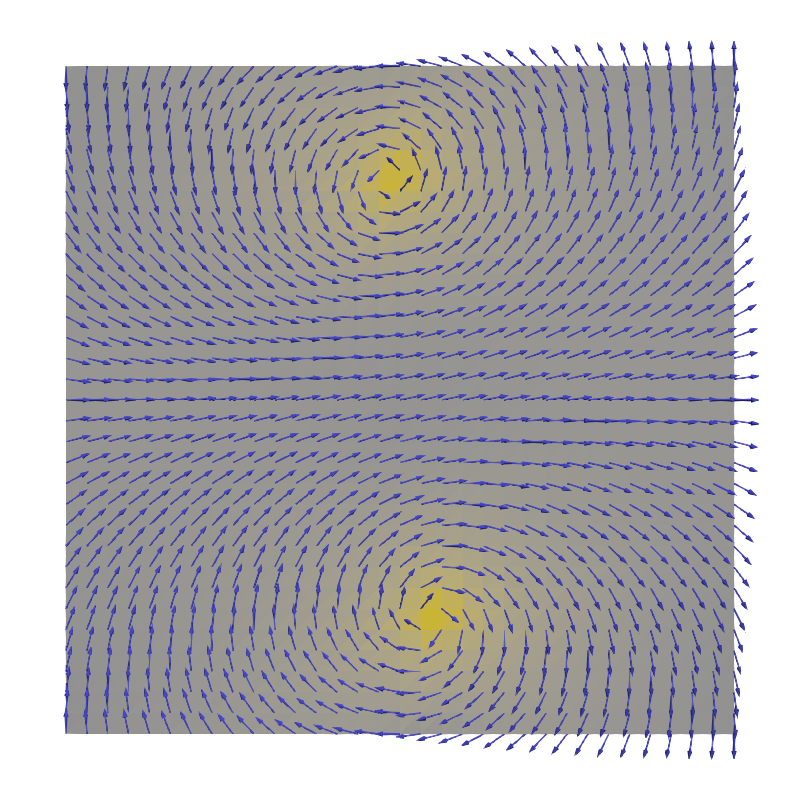}%
    \includegraphics[width=0.25\textwidth]{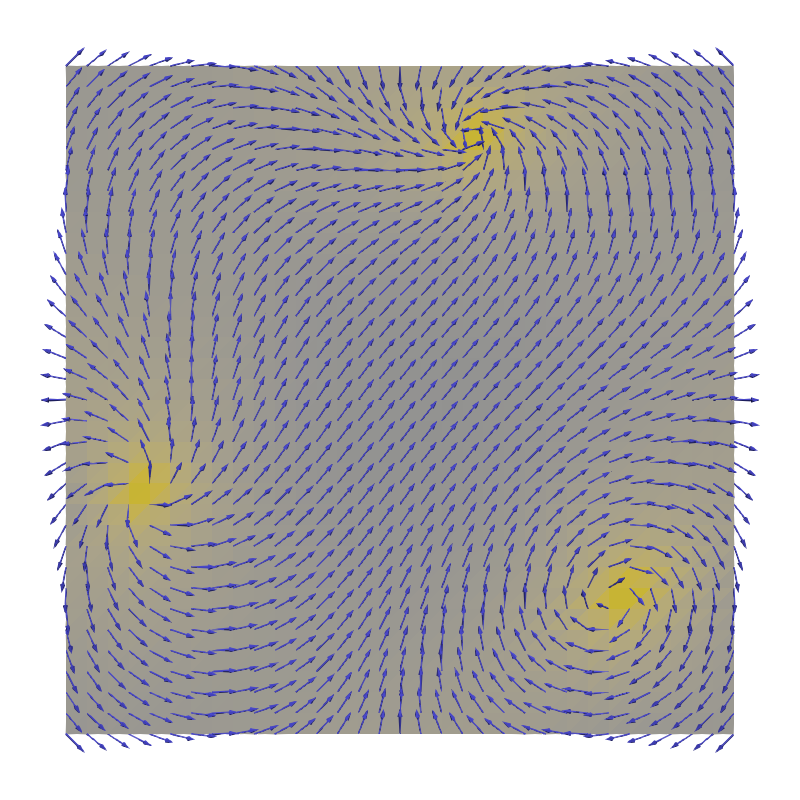}%
    \includegraphics[width=0.25\textwidth]{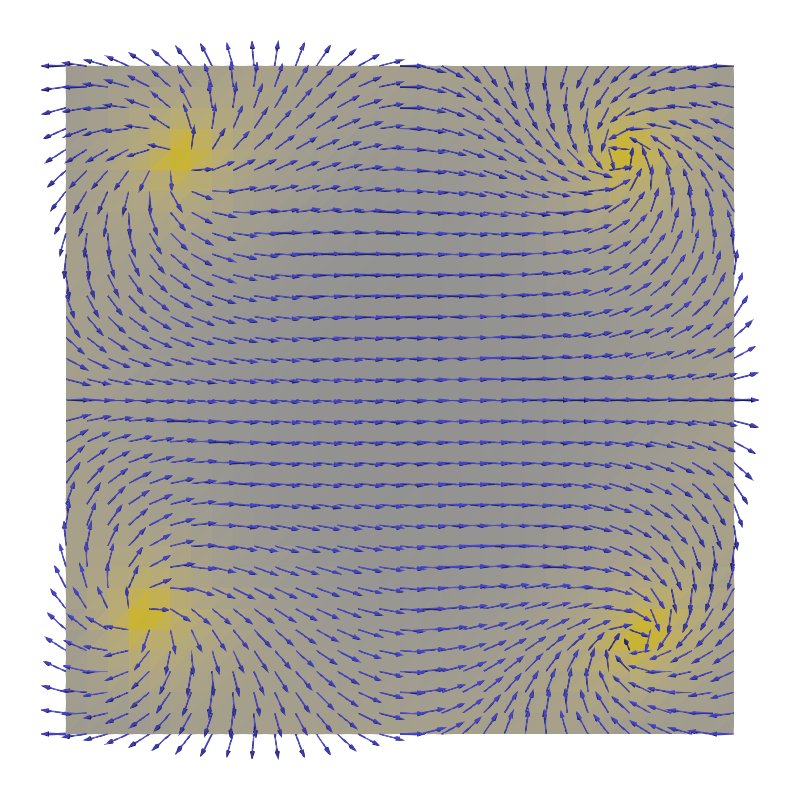}%
    \includegraphics[width=0.25\textwidth]{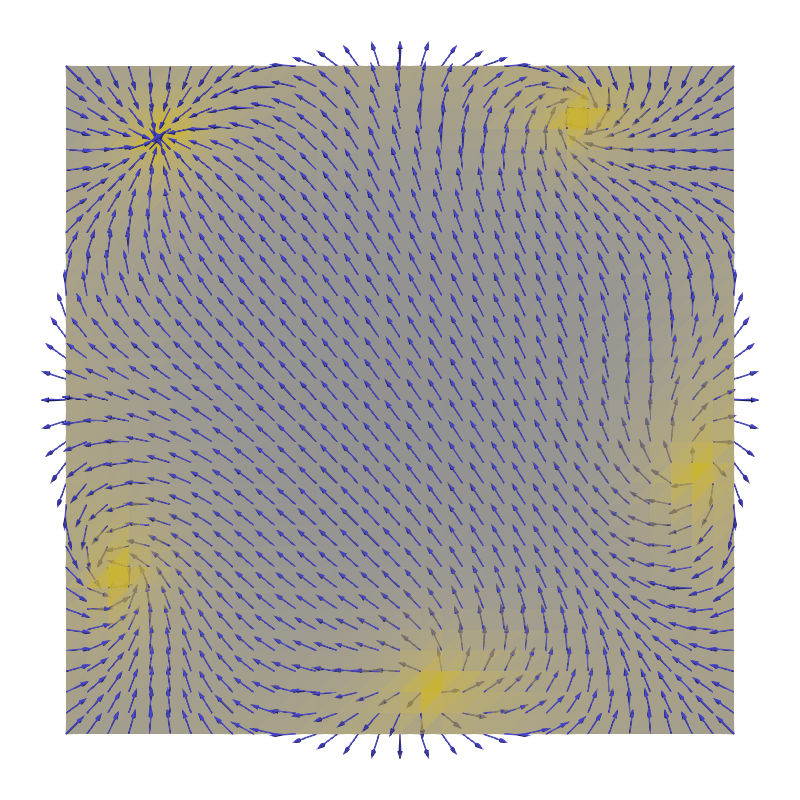}%
    \subcaption{Nonconforming discretization}
  \end{subfigure}

  \begin{subfigure}[b]{\textwidth}
    \includegraphics[width=0.25\textwidth]{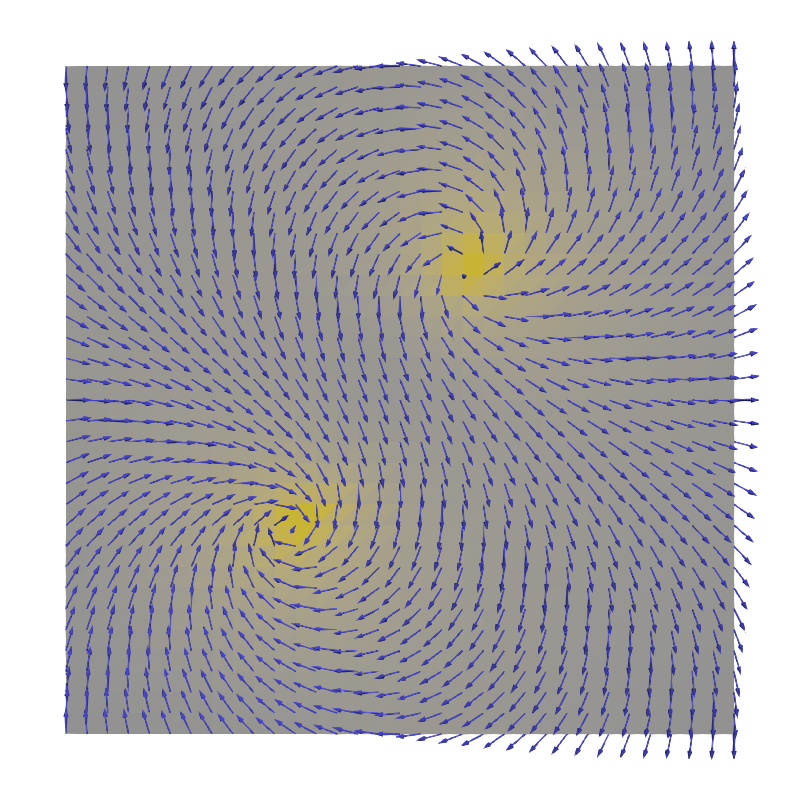}%
    \includegraphics[width=0.25\textwidth]{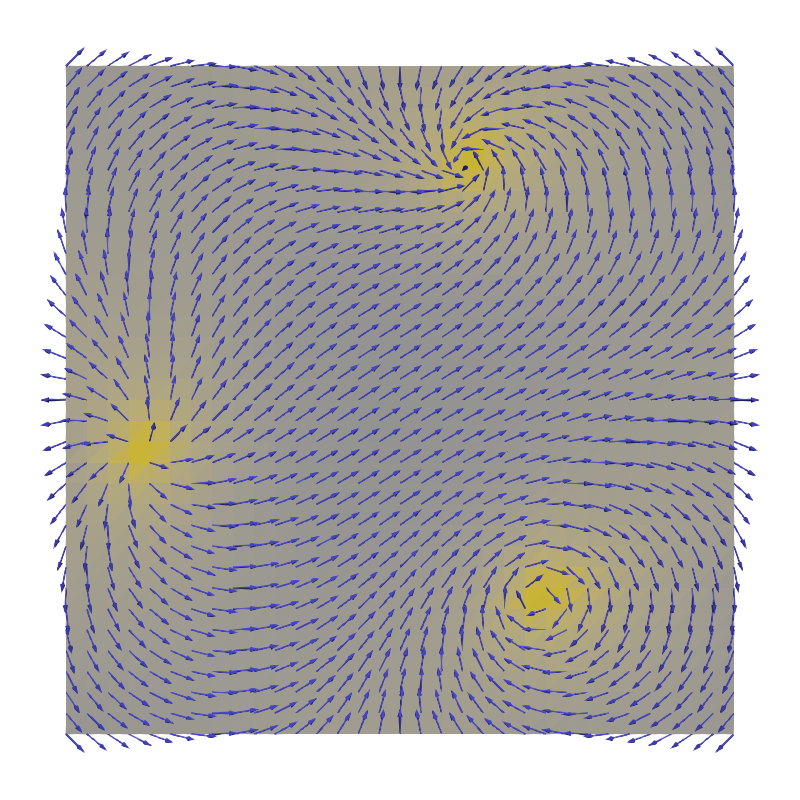}%
    \includegraphics[width=0.25\textwidth]{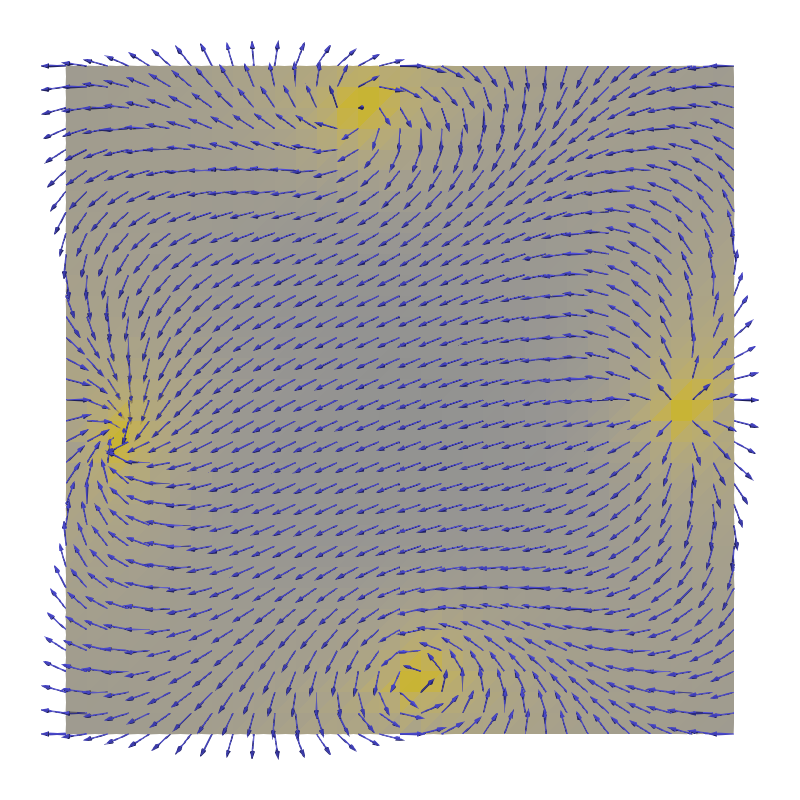}%
    \includegraphics[width=0.25\textwidth]{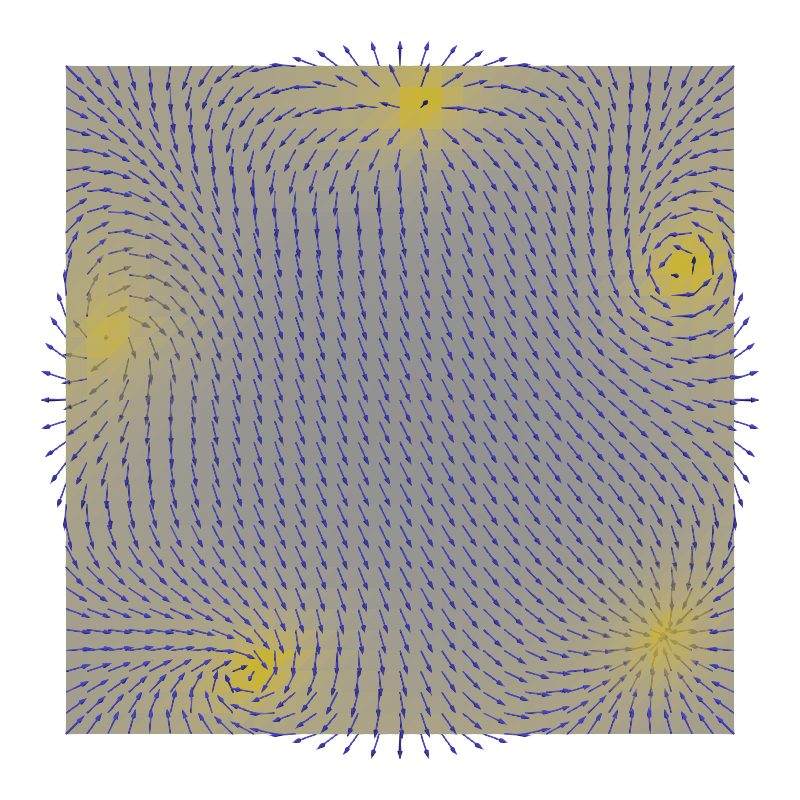}%
    \subcaption{Conforming discretization}
  \end{subfigure}
  \caption{Cut through midsurface $x_3=0$ for $\kappa = 2,\dots,5$ in Problem~\ref{prob:highsing},
   using the trust-region solver}
  \label{fig:highsing33_configurations}
\end{figure}

As the final test we measure iteration numbers and constraint violation for a harmonic map with several singularities (Problem~\ref{prob:highsing}).
Here, the start iterate $\bu_h^0$ is the function $\bu_D^\kappa$ defined in~\eqref{eq:singularity_kappa},
which has a singularity of degree~$\kappa$. As such a configuration
cannot stable for $\kappa \ge 2$ \cite{BreCorLie86}, we expect the solvers to converge to configurations
with $\kappa$ isolated singularities of degree~$1$. Figure~\ref{fig:highsing33_configurations}
shows that this does actually happen:  The final configuration consists of $\kappa$ singularities placed at roughly equal distances on a circle.  Unfortunately, one notices that
the final configurations differ considerably depending on which discretization is used.
Recalling from Chapter~\ref{ex:random} that the methods sometimes converge to configurations where the singularity was not at the
correct spot (in particular when using the conforming discretization), the configurations
shown in Figure~\ref{fig:highsing33_configurations} should not receive a lot of trust.
A deeper understanding of the behavior of discrete harmonic map models is needed before
this issue can be fully resolved.

Regarding the solver performance, the usual metrics are shown in Table~\ref{tab:highsing33-S}
for first-order discretizations.  There, one sees that numerically, this scenario does seem
to be more difficult than the others. In particular, for the trust-region methods the
iteration numbers are not bounded anymore, but increase slowly with increasing mesh size.
Unlike in earlier situations, more iterations are now needed for the conforming discretization
than for the nonconforming one.
The iteration numbers for the gradient-flow method are again much higher, but do not show a clear pattern as in earlier case.
The constraint violation remains in a reasonable range, and gets smaller with decreasing $\tau=4h$.

\begin{table}\footnotesize
\begin{center}
\begin{tabular}{c r | c c c| c  c| c c}
& &\multicolumn{3}{c|}{\textsf{Solve(GF)}/\textsf{Discr(N)}}   &\multicolumn{2}{c}{\textsf{Solve(TR)}/\textsf{Discr(N)}}  &\multicolumn{2}{c}{\textsf{Solve(TR)}/\textsf{Discr(C)}}  \\
  $r$ &$\abs{\Th}$ &\#Iter &$\delta_\mathrm{1}[\bu_h]$ & time [s] 
  &\#Iter &time [s] 
  &\#Iter & time [s] \\ \hline\hline
& & \multicolumn{7}{c}{degree $\kappa = 2$} \\
  1&  48        &   1      & \num{2.776e-17} &<1               & 1  & \timenum{  0.0308565} &  6 & \timenum{  0.112885}\\
  2&  384       &   238    & \num{6.920e-03} & \timenum{6}     & 8  & \timenum{  0.623017 } & 11 & \timenum{  0.664494}\\
  3&  3072      &   662    & \num{2.970e-03} & \timenum{27}    & 9  & \timenum{  2.08213  } & 20 & \timenum{  4.62172 }\\
  4&  24\,576   &   801    & \num{8.260e-04} & \timenum{296}   &21  & \timenum{ 25.4594   } & 17 & \timenum{ 23.2129  }\\
  5&  196\,608  &   3004   & \num{2.068e-04} & \timenum{13960} &19  & \timenum{164.47     } & 33 & \timenum{368.75    }\\
  \hline
& & \multicolumn{7}{c}{degree $\kappa = 3$} \\
  1&  48        &    1 & \num{7.864e-17} & <1               & 1 & \timenum{  0.708926} & 3 & \timenum{  0.0754735} \\
  2&  384       &  867 & \num{1.730e-02} & \timenum{22}     &12 & \timenum{  4.7088  } &19 & \timenum{  1.56621  } \\
  3&  3072      &  619 & \num{7.047e-03} & \timenum{26}     &11 & \timenum{  4.51344 } &16 & \timenum{  3.85301  } \\
  4&  24\,576   & 1310 & \num{2.040e-03} & \timenum{491}    &14 & \timenum{ 15.8183  } &29 & \timenum{ 39.0861   } \\
  5&  196\,608  & 3058 & \num{5.193e-04} & \timenum{14930}  &23 & \timenum{193.976   } &38 & \timenum{453.057    } \\
  \hline
& & \multicolumn{7}{c}{degree $\kappa = 4$}\\
  1&  48        &  23  & \num{4.919e-02} &<1                & 4  & \timenum{  0.0817272} & 7 & \timenum{  0.135057}  \\   
  2&  384       & 194  & \num{2.267e-02} & \timenum{5}      & 9  & \timenum{  1.37782  } &15 & \timenum{  1.91549 }  \\   
  3&  3072      &1198  & \num{1.129e-02} & \timenum{48}     &11  & \timenum{  4.02524  } &24 & \timenum{  5.65278 }  \\
  4&  24\,576   &1238  & \num{3.451e-03} & \timenum{483}    &27  & \timenum{ 33.6956   } &29 & \timenum{ 42.0638  }  \\  
  5&  196\,608  &2908  & \num{8.699e-04} & \timenum{15000}  &38  & \timenum{285.293    } &44 & \timenum{530.704   }  \\
  \hline
& & \multicolumn{7}{c}{degree $\kappa = 5$} \\
  1&  48        &   1  & \num{4.857e-17} & <1               & 1 & \timenum{0.71218}   & 3 & \timenum{0.0756177}   \\   
  2&  384       & 303  & \num{2.521e-02} & \timenum{7}      &13 & \timenum{2.16994}   &13 & \timenum{0.916537}    \\   
  3&  3072      & 828  & \num{1.535e-02} & \timenum{34}     &11 & \timenum{6.15935}   &25 & \timenum{6.80908}     \\
  4&  24\,576   &1197  & \num{4.775e-03} & \timenum{468}    &16 & \timenum{23.1106}   &33 & \timenum{43.0659}      \\  
  5&  196\,608  &3005  & \num{1.196e-03} & \timenum{19660}  &22 & \timenum{178.889}   &56 & \timenum{655.635}       \\
  \hline
\end{tabular}%
\end{center}
\caption{Iteration numbers, unit-length violation and wall-times for Problem~\ref{prob:highsing}}
\label{tab:highsing33-S}
\end{table}

As the number of required trust-region iterations has increased, the total wall-time difference
between the two solver algorithms is not as big anymore, at least for grid sizes in the range
that we could measure.  However, the difference is still considerable.

\medskip{}

\noindent \small{\textbf{Acknowledgements}
The authors gratefully acknowledge the support by the Deutsche Forschungs\-gemeinschaft
(Funder DOI: \url{https://dx.doi.org/10.13039/501100001659})
in the Research Unit~3013 \emph{Vector- and Tensor-Valued Surface PDEs} within the sub-projects \emph{TP3: Heterogeneous thin structures with prestrain} and \emph{TP4: Bending plates of nematic liquid crystal elastomers}.
}

\bibliographystyle{siam}
\bibliography{paper-bartels-boehnlein-palus-sander-harmonic-maps}

\end{document}